\documentclass[review]{elsarticle}
\usepackage[latin1]{inputenc}
\usepackage[T1]{fontenc}
\usepackage{amsmath}
\usepackage{amsfonts}
\usepackage{amssymb}
\usepackage{mathrsfs}
\usepackage{graphicx}          

\usepackage{tikz}
\usetikzlibrary{shapes,arrows}
\usetikzlibrary{decorations.text}

\journal{Chaos, Solitons and Fractals}

\newcommand{\SV}{SIRSi-Vaccine}
\newcommand{\cvd}{Covid-19}
\newcommand{\SP}{S\~ao Paulo}
\newcommand{\St}{Santos}
\newcommand{\Camp}{Campinas}

\begin{document}

\begin{frontmatter}

\title{\SV\ dynamical model for \cvd\ pandemic}

\author[Add1]{Cristiane M. Batistela}
\ead{cmbatistela@yahoo.com.br}

\author[Add2]{Diego P. F. Correa}
\ead{diego.ferruzzo@ufabc.edu.br}

\author[Add3]{\'Atila M Bueno}
\ead{atila.bueno@unesp.br}

\author[Add1]{Jos\'e R. C. Piqueira\corref{mycorrespondingauthor}}
\cortext[mycorrespondingauthor]{Corresponding author}
\ead{piqueira@lac.usp.br}

\address[Add1]{ S\~ao Paulo University, Polytechnic School - EPUSP, S\~ao Paulo, SP, Brazil.}
\address[Add2]{Federal University of ABC - UFABC, S\~ao Bernardo do Campo, SP, Brazil}
\address[Add3]{S\~ao Paulo State University - UNESP, Sorocaba, SP, Brazil}

\begin{abstract}
The Severe  Acute Respiratory Syndrome Corona Virus 2 (SARS-CoV-2), or \cvd, burst 
into a pandemic in the beginning of 2020. An unprecedented worldwide effort involving 
academic institutions, regulatory agencies and industry is facing the challenges 
imposed by the rapidly spreading disease. Emergency use authorization for vaccines 
were granted in the beginning of December 2020 in Europe and nine days later in the United States. 
The urge for vaccination started a race, forcing governs and health care agencies to take decisions on the fly  
regarding the vaccination strategy and logistics. So far, the vaccination strategies and  
non-pharmaceutical interventions, such as social distancing and the use of face masks, are the only efficient 
actions to stop the pandemic.  In this context, it is of fundamental importance to understand 
the dynamical behavior of the \cvd\ spread along with possible vaccination strategies.
In this work a Susceptible - Infected - Removed - Sick with vaccination (\SV) model is proposed. 
In addtion, the \SV\ model also accounts for unreported, or asymptomatic, cases and the possibility 
of temporary immunity, either after infection or vaccination.  Disease free and endemic equilibrium 
points existence conditions are determined in the ($\omega\times\theta$) vaccine-effort and social distancing 
parameter space. The model is adjusted to the data from \SP, \St\ and \Camp, three major cities in the 
State of S\~ao Paulo, Brazil. 
\end{abstract}

\end{frontmatter}

\section{Introduction}
\label{sec:Intro}
In 30 January 2020 the World Health Organization (WHO) declared a Public Health Emergency 
of International Concern (PHEIC). Four weeks before, the Wuhan Municipal Health Commission 
had reported a cluster of pneumonia cases, whose causative agent was soon after discovered to 
be the Severe  Acute Respiratory Syndrome Corona Virus 2 (SARS-CoV-2), or \cvd. Despite the 
WHO alert,  the \cvd\ became pandemic weeks later, in March 2020  
\cite{Maxmen_2021,WHOTimelineCovid19}. 

It is of great concern the fact that worldwide governs and Health Care Agencies were unable to stop 
the Covid-19 pandemic in its beginning. As a result, the rapidly spreading disease is challenging 
health care angencies, academic institutions and industry to develop and deploy efficient drug 
treatments and vaccination. 

A number of drugs have been tried in patients with \cvd\ disease. Unfortunately, studies 
\cite{doi:10.1056/NEJMoa2023184,JOHNSTON2021100773,Tomlinson:2021vp,Bush:2020vu,Bignardi_2021} 
found that the drugs had little or no effect on the overall mortality, initiation of ventilation, duration of hospital stay  
or viral clearance even for patients with chronic use of some of the tried drugs.  

On the other hand, non-pharmaceutical interventions (NPI), such as physical distancing, use of face masks and 
eye protection, reduce the virus transmission \cite{CHU20201973,DeSouzaSantos2021}. 
Worldwide populations have been compliant to NPI, however,  in some countries, there is a controversy 
over its effectiveness \cite{ijerph17186655}. In Brazil, the closure of non-essential activities lasted only for 
a short time, and the cities lifted NPI in an uncoordinated manner \cite{DeSouzaSantos2021}. As a result,  
in the beginning of November 2020,  the downward trend reversed and the number of \cvd\ cases started 
to raise in a second wave of infection \cite{PainelCoronavirusMarco2021}.

The development of vaccines usually takes many years, even decades, and is also very expensive.  
The fast development of \cvd\ vaccines was only possible since researchers had been working for 
years on vaccines for similar viruses, such as, SARS (Severe Acute Respiratory Syndrome) and MERS 
(Middle East Respiratory Syndrome). In addition, the experience gained with the Ebola vaccine 
showed that the development of new vaccines can be accelerated with a worldwide effort, including 
academic institutions, industry and health care regulatory agencies, without compromising safety
\cite{Ball:2021uk,Wolf:2020vo,FinneganVaccinesToday2021}.

There are currently 78 vaccines in clinical trial on humans, and 
other 77 are under investigation in animals. Emergency use authorization for vaccines 
were granted in the beginning of December 2020 by Britain's National Health Service, and 
later on  by the U.S. Food and Drug Administration (FDA). In january 2021 the Brazilian national health 
surveillance agency (Anvisa) also granted emergency use for vaccines.  
Up to now, in several countries, 6 vaccines have been granted emergency use, 
and 7 approved for full use   \cite{NYT-VaccineTracker,AnvisaAprovaVacina}. 

The \cvd\ pandemic, and the consequent urge for vaccination, started a race, with countries trying 
to vaccinate populations as fast as possible. In \cite{doi:10.1056/NEJMoa2101765} the immunization 
results of a nationwide mass vaccination in Israel is reported, strengthening the expectation 
of \cvd\ pandemic effects mitigation. However, at the current pace, and for most countries,  
a successful vaccination strategy consists of a long journey ahead  \cite{OurWorldInData}. 

Logistics plays an important role in vaccination strategies deployment, since the efficacy\footnote{The vaccine efficacy 
is the percentage of reduction of a disease in the vaccinated group in a clinical trial.}  of 
the vaccines currently in development ranges from 50.38\% to 95.0\% \cite{NYT-VaccineTracker}, 
and that different vaccines need different vaccination strategies, given the number of doses, the 
number of weeks apart each dose, specific storage facility needs, and so on. 

Besides, Health Care Agencies must take decisions considering the \cvd\ new variants, 
vaccination strategy with varying NPI adherence among the population, or consider to delay 
the second dose in order to speed up the one dose vaccinated population  
\cite{10.1001/jama.2021.1114,doi:10.1056/NEJMclde2101987}. In addition, important factors 
are vaccine efficacy and population vaccination cover. Vaccines have to have at least $70\%$ efficacy 
with population cover of $75\%$ to prevent a pandemic without any other intervention. 
To stop the ongoing pandemic the needed vaccines efficacy raises to at least $80\%$, for $75\%$ 
population cover  \cite{BARTSCH2020493}. 

Given such complex context it is important to understand the dynamical behavior of the 
\cvd\ pandemic, considering possible vaccination strategies in order to mitigate of even 
to extinguish the ongoing pandemic.    

Many scientific works have been developed since the beginning of the pandemic, 
aiming to stop, or mitigate, the pandemic effects, consisting of control systems strategies, 
and deep learning modeling. Other works, consisting of mathematical modeling, analysis, 
and numerical simulations, in order to predict the \cvd\ contagious behavior, have also been 
published.  

In \cite{ARUNKUMAR2021107161} time-series models - Auto-Regressive Integrated Moving Average (ARIMA) 
and Seasonal Auto-Regressive Integrated Moving Average (SARIMA) are developed for the COVID-19 pandemic, 
with data for the countries where $70\%$ - $80\%$ of global cumulative cases are located. The results show 
exponential growth of confirmed cases and deaths. In \cite{FIACCHINI2021} a two parameter model is used to 
predict the pandemic evolution in France. 

A Susceptible - Infected - Removed model is used to predict the \cvd\ pandemic in Kuwait \cite{ALENEZI20213161}. 
A SIR-type model with non-constant population is developed In \cite{MUNOZFERNANDEZ2021110682}.  
The model is calibrated to the rates of infection and deaths for Italy, Spain and the USA. 

In \cite{NAKAMURA2021110667} the death toll due to the \cvd\ pandemic is modeled via sigmoid curve,  
predicting mortality, and estimating practical interest values such as the number of new cases per 
infection, i.e., the basic reproduction number. 

In \cite{WEI2021} a mathematical model of \cvd\ is developed. A control system is proposed 
aiming to reduce the number of the susceptible individuals in a population applying NPI  
consisting of social distancing.  

A number of works consider deep learning to predict the number of \cvd\ cases. 
In \cite{GUPTA2021107039} a prediction model of confirmed and death cases of \cvd\ is developed, 
based on a deep learning algorithm with two long short-term memory layers, considering 
the available data for \cvd\ in India.  In \cite{Blyuss:2021wp}, multilayer deep learning models with perceptron, 
random forest and Long short-term memory, are trained and used to forecast the \cvd\ 
pandemic in Iran.   

In \cite{BATISTELA2021110388} a Susceptible - Infected - Removed - Sick (SIRSi) model is proposed,
modeling unreported or asymptomatic cases, and considering the effects of temporary acquired immunity 
on the \cvd\ spread. 

In this work, the model proposed in \cite{BATISTELA2021110388} is modified, and a vaccination 
strategy is included, consisting of a vaccination rate of susceptible individuals. 
The proposed \SV\ model is composed of four compartments, namely, Susceptible, Infected, 
Sick and Recovered. The focus is to assess the effects of vaccination associated with social 
distancing as a NPI in the \cvd\ spread dynamics.  

The proposed \SV\ model presents both disease free and endemic equilibrium points. 
The equilibrium points stabilities are studied both analytically and numerically, and, as a result,  
the found equilibrium existence condition relates the effects of the social distancing  
to the vaccination rate. 

In addition, the \SV\ model is numerically adjusted to the epidemic situation in 
\SP, \St\ and \Camp, three major cities in the State of \SP\ connecting the 
shores with the interior of the S\~ao Paulo state. The numerical results show that 
in order to stop the \cvd\ spread different strategies may be necessary for 
different cities, concerning social distancing NPI and vaccination rates. 

The paper is organized as follows; Section \ref{sec:MathModel} presents the 
mathematical model. In Section \ref{sec:EqPoints} the equilibrium points are 
determined and their stability are investigated. The parametric fitting 
is presented in Section \ref{sec:ParamFit}.  In Section \ref{sec:SimResults} 
the simulations results for \SP, \St\ and \Camp\ are presented. 
The final remarks are shown in Section \ref{sec:Conclusions}.

\section{SIRSi and \SV\ mathematical models}
\label{sec:MathModel}
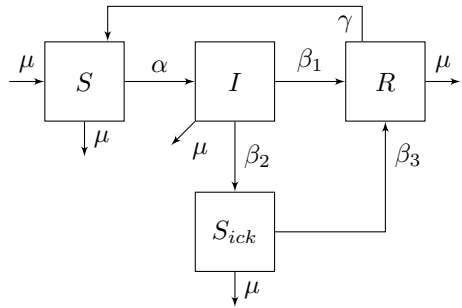
\begin{figure}[!htb]
\tikzstyle{block} = [draw, fill=none, rectangle, minimum height=3em, minimum width=3em, node distance=2cm]
\tikzstyle{sum} = [draw, fill=white, circle, node distance=auto]
\tikzstyle{input} = [coordinate]
\tikzstyle{output} = [coordinate]
\tikzstyle{pinstyle} = [pin edge={to-,thin,black}]
\begin{center}
\begin{tikzpicture}[auto, node distance=2cm,>=latex']
    \node [input] (In){};
    \node [block, right of=In, node distance=1cm] (S) {$S$};
    \node [block, right of=S] (I) {$I$};
	\node [output, below of=S, node distance=1cm] (d) {};    
    \node [block, right of=I] (R) {$R$};
    \node [output, right of=R, node distance=1cm] (Im) {};
    \node [block, below of=I] (Sick) {$S_{ick}$};
    \node [output, below of=Sick, node distance=1cm] (D) {};    
    \node [output, below left of=I, node distance=1.2cm] (D2) {};
    \draw [draw,->] (In)-- node {$\mu$} (S);
    \draw [draw,->] (S)-- node {$\mu$} (d);
    \draw [draw,->] (S) -- node {$\alpha$} (I);
	\draw [draw,->] (I) -- node {$\beta_1$} (R);
    \draw [draw,->] (R) -- node {$\mu$} (Im);
    \draw [draw,->] (R.120) |- node[at start, anchor =south east]{$\gamma$} (2,1.0) -| (S.60);
    \draw [draw,->] (I) -- node {$\beta_2$} (Sick);
    \draw [draw,->] (I) -- node {$\mu$} (D2);    
    \draw [draw,->] (Sick) -- node {$\mu$} (D);    
    \draw [draw,->] (Sick) -| node[near end, anchor=south west]{$\beta_3$} (R);
\end{tikzpicture}
\end{center}
\caption{SIRSi model with constant population.}
\label{fig:sirsi-model}
\end{figure}

A constant population model with four compartments is proposed, see  Fig. \ref{fig:sirsi-model}. 
The compartments are defined as:  Susceptible $(S)$, Infected $(I)$, 
Infected symptomatic positive tested $(S_{ick})$ and Recovered $(R)$. 

In Fig \ref{fig:sirsi-model}, $I$ is the infectious population compartment representing the population 
in the incubation stage, i.e., prior to the onset of symptoms. Transmission during this period 
has been reported in \cite{doi:10.1056/NEJMc2001737,Zhang2020,Tian2020}. 

The infected population can be asymptomatic or symptomatic, and according to \cite{Ferguson2020}   
the incubation period is assumed to be of 5.1 days. Infectiousness occur $12$ hours prior to the symptoms onset 
for the symptomatic. On the other hand, for those who are asymptomatic, the infectiousness is assumed 
to occur $4.6$ days after infection. The average time between infection and infectiousness is 6.5 days. 

Those who are asymptomatic or do not develop severe symptoms, i.e., neither tested nor  
documented cases, are moved to the  $R$ (Recovered) compartment after a $1/\beta_{1}$ 
period \cite{BATISTELA2021110388}.

The incubation period for those who are sympotmatic is $1/\beta_2$. Once the infected individual is tested 
positive and the case is documented, it is moved to the $S_{ick}$ compartment, which consists of 
those patients with severe symptoms seeking medical attention. 

After a period $1/\beta_3$ the $S_{ick}$ population that recovers is moved to the compartment $R$. 

The possibility of temporary immunity for recovered individuals, as addressed 
in \cite{mbs:/content/journal/jgv/10.1099/jgv.0.001439,Callow:1990vc,Mo2006,BATISTELA2021110388},  
is also considered. As it can be seen in Fig. \ref{fig:sirsi-model}, the recovered $R$ population 
becomes susceptible $S$ again at a rate $\gamma$. 

The model is developed considering the population size constraint of Eq. \ref{eq:cond-const-pop}. 
In addition, the population growth and death rates $\mu$, including deaths due to \cvd, are  
similar in the State of São Paulo \cite{SEADESocDist}, therefore, in this work, they are considered to be equal. 

\begin{equation}
N=S(t)+I(t)+S_{ick}(t)+R(t) 
\label{eq:cond-const-pop}
\end{equation}

From the foregoing considerations the  SIRSi mathematical model is given by Eq. \ref{eq:constant-pop-model}.

\begin{equation}
\begin{split}
\frac{d}{dt}S(t)&=\mu N-\frac{\alpha(1-\theta)}{N}S(t)I(t)-\mu S(t)+\gamma R(t);\\
\frac{d}{dt}I(t)&=\frac{\alpha(1-\theta)}{N}S(t)I(t)-(\beta_1+\beta_2)I(t)-\mu I(t);\\
\frac{d}{dt}S_{ick}(t)&=\beta_2I(t)-\beta_3S_{ick}(t)-\mu S_{ick}(t);\\
\frac{d}{dt}R(t)&=\beta_1I(t)+\beta_3S_{ick}(t)-\gamma R(t) -\mu R(t);
\end{split}
\label{eq:constant-pop-model}
\end{equation}

In Eq. \ref{eq:constant-pop-model}, the effect of social distancing measure, which is 
a form of NPI, is introduced by the parameter $\theta$, $0<\theta<1$.  
Where $\theta=0$ means that no social distancing measure is under consideration,   
and $\theta=1$ means a complete lockdown.   

A normalization with respect to the population size $N$ is performed in Eq. \ref{eq:cond-const-pop}, 
resulting in Eq. \ref{eq:cond-const-pop-NORM1},

\begin{equation}
1=s(t)+i(t)+s_{ick}(t)+r(t) 
\label{eq:cond-const-pop-NORM1}
\end{equation}
 
\noindent where $s(t)=S(t)/N$, $i(t)=I(t)/N$, $s_{ick}(t)=S_{ick}(t)/N$ and $r(t)=R(t)/N$. 
Replacing Eq. \ref{eq:cond-const-pop-NORM1} into the system of Eq. \ref{eq:constant-pop-model},  
results Eq. \ref{eq:constant-perc-pop-model}.

\begin{equation}
\begin{split}
\frac{d}{dt}s(t)&=\mu -\alpha(1-\theta)si-\mu s+\gamma r;\\
\frac{d}{dt}i(t)&=\alpha(1-\theta)si-(\beta_1+\beta_2)i-\mu i;\\
\frac{d}{dt}s_{ick}(t)&=\beta_2i-\beta_3s_{ick}-\mu s_{ick};\\
\frac{d}{dt}r(t)&=\beta_1i+\beta_3s_{ick}-\gamma r -\mu r;
\end{split}
\label{eq:constant-perc-pop-model}
\end{equation}

In addition, from Eq. \ref{eq:cond-const-pop-NORM1}, the recovered  compartment $r(t)$ 
can be written as a linear combination of the other compartments (or state variables), 
as shown in Eq. \ref{eq:cond-const-pop-NORM2}, 

\begin{equation}
r(t)=1-s(t)-i(t)-s_{ick}(t),  
\label{eq:cond-const-pop-NORM2}
\end{equation}

\noindent consequently, the solutions $\Omega=\{(s(t),i(t),s_{ick}(t))\in\mathbb{R}^3_+\,|\,s(t)+i(t)+s_{ick}(t)\leq 1\}$ 
of the system of Eq. \ref{eq:constant-perc-pop-model} can be studied in the system 
of Eq. \ref{eq:constant-perc-pop-reduced-model4}, 

\begin{equation}
\begin{split}
\frac{d}{dt}s(t)&=\mu+\gamma -\alpha(1-\theta)si-(\mu+\gamma) s -\gamma i -\gamma s_{ick}\\
\frac{d}{dt}i(t)&=\alpha(1-\theta)si-(\beta_1+\beta_2+\mu)i\\
\frac{d}{dt}s_{ick}(t)&=\beta_2i-(\beta_3+\mu)s_{ick}
\end{split}
\label{eq:constant-perc-pop-reduced-model4}
\end{equation}

\noindent with reproduction basic number,

\begin{equation}
R_0=\dfrac{\alpha(1-\theta)}{\beta_1+\beta_2+\mu}.
\label{eq:R0-model1}
\end{equation}

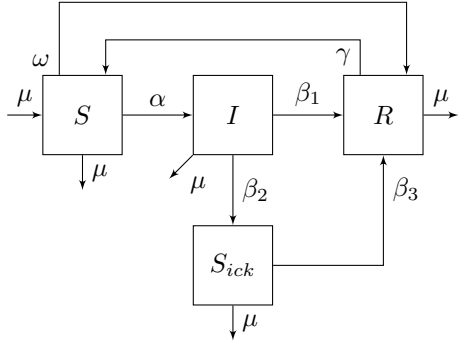
\begin{figure}[!htb]
\tikzstyle{block} = [draw, fill=none, rectangle, minimum height=3em, minimum width=3em, node distance=2cm]
\tikzstyle{sum} = [draw, fill=white, circle, node distance=auto]
\tikzstyle{input} = [coordinate]
\tikzstyle{output} = [coordinate]
\tikzstyle{pinstyle} = [pin edge={to-,thin,black}]
\begin{center}
\begin{tikzpicture}[auto, node distance=2cm,>=latex']
    \node [input] (In){};
    \node [block, right of=In, node distance=1cm] (S) {$S$};
    \node [block, right of=S] (I) {$I$};
	\node [output, below of=S, node distance=1cm] (d) {};    
    \node [block, right of=I] (R) {$R$};
    \node [output, right of=R, node distance=1cm] (Im) {};
    \node [block, below of=I] (Sick) {$S_{ick}$};
    \node [output, below of=Sick, node distance=1cm] (D) {};    
     \node [output, below left of=I, node distance=1.2cm] (D2) {};
    \draw [draw,->] (In)-- node {$\mu$} (S);
    \draw [draw,->] (S)-- node {$\mu$} (d);
    \draw [draw,->] (S) -- node {$\alpha$} (I);
	\draw [draw,->] (I) -- node {$\beta_1$} (R);
    \draw [draw,->] (R) -- node {$\mu$} (Im);
    \draw [draw,->] (R.120) |- node[at start, anchor =south east]{$\gamma$} (2,1.0) -| (S.60);
    \draw [draw,->] (S.120) |- node[at start, anchor =south east]{$\omega$} (2,1.5) -| (R.60);
    \draw [draw,->] (I) -- node {$\beta_2$} (Sick);
    \draw [draw,->] (Sick) -- node {$\mu$} (D);    
    \draw [draw,->] (Sick) -| node[near end, anchor=south west]{$\beta_3$} (R);
    \draw [draw,->] (I) -- node {$\mu$} (D2);    
\end{tikzpicture}
\end{center}
\caption{SIRSi model with vaccination intervention.}
\label{fig:sirsi-model-vacc}
\end{figure}

A vaccination intervention is introduced consisting of susceptible individuals 
vaccinated at a rate \(\omega\). The vaccinated individuals are considered to become 
instantaneously immune after vaccination, being moved to the recovered compartment. 
The proposed \SV\ model is shown in Fig. \ref{fig:sirsi-model-vacc}.

Including the vaccination intervention in Eqs. \ref{eq:constant-perc-pop-model} and  
\ref{eq:constant-perc-pop-reduced-model4} results in Eqs. \ref{eq:constant-perc-pop-reduced-model2} 
and \ref{eq:constant-perc-pop-reduced-model3}, respectively. 
	
\begin{equation}
\begin{split}
\frac{d}{dt}s(t)&=\mu -\alpha(1-\theta)si-\mu s+\gamma r-\omega s\\
\frac{d}{dt}i(t)&=\alpha(1-\theta)si-(\beta_1+\beta_2)i-\mu i\\
\frac{d}{dt}s_{ick}(t)&=\beta_2i-\beta_3s_{ick}-\mu s_{ick}\\
\frac{dr}{dt}r(t)&=\beta_1i+\beta_3s_{ick}-\gamma r +\omega s-\mu r
\end{split}
\label{eq:constant-perc-pop-reduced-model2}
\end{equation}

\begin{equation}
\begin{split}
\frac{d}{dt}s(t)&=\mu+\gamma -\alpha(1-\theta)si-(\mu+\gamma+\omega) s -\gamma i -\gamma s_{ick}\\
\frac{d}{dt}i(t)&=\alpha(1-\theta)si-(\beta_1+\beta_2+\mu)i\\
\frac{d}{dt}s_{ick}(t)&=\beta_2i-(\beta_3+\mu)s_{ick}
\end{split}
\label{eq:constant-perc-pop-reduced-model3}
\end{equation}

\section{Disease-free and endemic equilibrium points}
\label{sec:EqPoints}
In this section the equilibrium points for the mathematical models 
in Eqs. \ref{eq:constant-perc-pop-reduced-model4} and \ref{eq:constant-perc-pop-reduced-model3} 
are determined and their stabilities are analysed. 

\subsection{SIRSi equilibrium points}
\label{sec:EqPointsSIRSi}
The system of Eq. \ref{eq:constant-perc-pop-reduced-model4} has two equilibrium points,   
one is consistent with a disease-free situation, and the other representing an endemic equilibrium.

The disease-free equilibrium point $P_{df1}$ is given by Eq. \ref{eq:PtoEq-P1}.  

\begin{equation}
P_{df1}=\left[\begin{array}{c}
s^{\ast}\\
i^{\ast}\\
s^{\ast}_{ick}
\end{array}\right]
=\left[\begin{array}{c}
1\\
0\\
0
\end{array}\right]
\label{eq:PtoEq-P1}
\end{equation}

On the other hand, the endemic equilibrium point $P_{e1}$ is given by Eq. \ref{eq:PtoEq-P2} 

\begin{equation}
P_{e1}=
\left[\begin{array}{c}
s^{\ast}\\
i^{\ast}\\
s^{\ast}_{ick}
\end{array}\right]= 
\left[
\begin{array}{c}
\frac{\beta_{1} + \beta_{2} + \mu}{\alpha \left(1-\theta \right)}\\
(\beta_3+\mu)\varphi\\
\beta_2\varphi
\end{array}
\right], 
\label{eq:PtoEq-P2}
\end{equation}
 
\noindent where 

\begin{equation}
\varphi=
\frac{\left(\gamma + \mu\right) \left(\alpha(1- \theta) - (\beta_{1} + \beta_{2} + \mu)\right)}
{\alpha \left(1-\theta \right) \left( (\beta_{1} + \beta_{2} + \mu)(\beta_3 +\mu) + (\beta_{2}  + \beta_{3} + \mu)\gamma  \right)}.
\label{eq:varphi}
\end{equation}
 
In order to analyze the local stability of the point, the Jacobian matrix $J$  of the 
system in Eq. \ref{eq:constant-perc-pop-reduced-model4} is shown  in Eq. \ref{eq:JacobianP1}. 

\begin{equation}
J = \left[
\begin{array}{@{}c@{}c@{}c@{}}
-\alpha(1-\theta)i^{\ast} -(\mu+\gamma) & -\alpha(1-\theta)s^{\ast} -\gamma & -\gamma \\
\alpha(1-\theta)i^{\ast} & \alpha(1-\theta)s^{\ast} - (\beta_1+\beta_2+\mu) & 0 \\
0 & \beta_2 & -(\beta_3 +\mu) 
\end{array}
\right]
\label{eq:JacobianP1}
\end{equation}

The eigenvalues associated with $P_{df1}$, the disease free equilibrium point, are 
given by Eq. \ref{eq:eig_free_disease_eq_mode}. 

\begin{equation}
\begin{split}
\lambda_1 &= -(\mu + \gamma)\\
\lambda_2 &=  \alpha(1-\theta)-(\beta_1+\beta_2+\mu)\\
\lambda_3 &= - (\beta_3 + \mu)
\end{split}
\label{eq:eig_free_disease_eq_mode}
\end{equation}

The eigenvalue $\lambda_2$ 
shows the possibility of a bifurcation since it can change sign depending on the 
parameters values. Therefore, the disease free equilibrium point $P_{df1}$ 
is stable if 

\begin{equation}
\alpha(1-\theta)<\beta_1+\beta_2+\mu.  
\label{eq:EstbilP1}
\end{equation}

Comparing Eqs. \ref{eq:R0-model1} and Eq. \ref{eq:EstbilP1} it can be noticed that 
if $R_0<1$ then $P_{df1}$ is locally asymptotically stable. On the other hand, if $R_0>1$, then 
$P_{df1}$ is unstable. 

Considering the endemic equilibrium point $P_{e1}$, Eqs. \ref{eq:PtoEq-P2} and \ref{eq:varphi}, 
and the Jacobian matrix, Eq. \ref{eq:JacobianP1}, the characteristic polynomial of $P_{e1}$ 
is given by Eq. \ref{eq:CharPolP2}, where the coefficients are given by Eq. \ref{eq:CharPolP2-1}. 

\begin{equation}
\mathcal{P}_{e1}(\lambda)=\lambda^{3}+a_{1}\lambda^2+a_{2}\lambda+a_{3}, 
\label{eq:CharPolP2}
\end{equation}

\begin{equation}
\begin{split}
a_{1}&=\varphi
\alpha(1- \theta)
\left( \beta_{3}  + \mu \right)+
\beta_{3} + \gamma + 2 \mu \\ 
a_{2}&=\left(\beta_{3} + \mu\right)
\left[\gamma + \mu+
\varphi\alpha(1-\theta)
\left(\beta_{1} + \beta_{2} + \beta_{3} + \gamma  + 2 \mu \right)\right]\\ 
a_{3}&=\varphi
\alpha\left(1-\theta\right) 
\left(\beta_{3} + \mu\right)
\left(
\begin{array}{@{}c@{}}
\beta_{1} \beta_{3} + \beta_{1} \mu + \beta_{2} \beta_{3} + \beta_{2} \gamma\\ 
+\beta_{2} \mu + \beta_{3} \gamma + \beta_{3} \mu + \gamma \mu + \mu^{2}
\end{array}
\right)
\end{split}
\label{eq:CharPolP2-1}
\end{equation}

The coefficients $a_{1}$,  $a_{2}$ and $a_{3}$, in Eq. \ref{eq:CharPolP2-1}, are positive real 
numbers if $\varphi>0$ (see Eq. \ref{eq:varphi}). In this case, in order to determine whether the 
characteristic polynomial in Eq. \ref{eq:CharPolP2} has unstable or stable roots the Routh-Hurwitz 
stability criterion \cite{Golnaraghi:2017vb} is used.  

Consequently, the characteristic polynomial in Eq. \ref{eq:CharPolP2} is stable if  
expression \ref{eq:RH1} holds.  
 
\begin{equation}
b_{1}=a_{2}-\frac{a_{3}}{a_{1}}>0
\label{eq:RH1}
\end{equation}

The proof that expression \ref{eq:RH1} holds true demands some mathematical reasoning, 
which results in Eq. \ref{eq:RH2}. 
 
\begin{equation}
\begin{split}
b_1&=\frac{1}{a_{1}}\left(\beta_{3} + \mu\right)
\left[\varphi^{2} c_{1}^{2}c_{2}+\varphi c_{1}c_{3}+c_{4}\right], 
\end{split}
\label{eq:RH2}
\end{equation}

\noindent where 

\begin{displaymath}
\begin{split}
c_{1}&=\alpha(1 -\theta),\\
c_{2}&=(\beta_{1}+\beta_2) (\beta_{3} + \mu)+ \beta_{3}^{2}+ 
\beta_{3} \gamma  + 3 \beta_{3} \mu  + \gamma \mu  + 2 \mu^{2},\\
c_{3}&=\beta_{1} \gamma + \beta_{1} \mu + \beta_{2} \mu + \beta_{3}^{2} 
+ 2 \beta_{3} \gamma + 4 \beta_{3} \mu + \gamma^{2} + 4 \gamma \mu + 4 \mu^{2},\\
c_{4}&=\beta_{3} \gamma + \beta_{3} \mu + \gamma^{2} + 3 \gamma \mu + 2 \mu^{2}.
\end{split}
\end{displaymath}

The expression inside brackets in Eq. \ref{eq:RH2} is positive if $\varphi>0$. Consequently, 
for $\varphi>0$ the expression \ref{eq:RH1} holds true, and  $P_{3}$ is asymptotically stable. 
In addition, as it can be seen in Eq. \ref{eq:varphi}, $\varphi>0$ implies that $R_{0}>1$.  

\subsection{\SV\ equilibrium points}
\label{sec:EqPointsSIRSiVaccine}
Similarly to what was found in Section \ref{sec:EqPointsSIRSiVaccine}, 
the \SV\ mathematical model of Eq. \ref{eq:constant-perc-pop-reduced-model3}  
has two equilibrium points, the first one is consistent with a disease-free situation, 
and the other with an endemic equilibrium.  

The disease free equilibrium point of Eq. \ref{eq:constant-perc-pop-reduced-model3} 
is shown in Eq. \ref{eq:free-disease-equilibrium}, for $s^{\ast}>0$. 

\begin{equation}
P_{df2}=
\left[
\begin{array}{c}
s^{\ast}\\
i^{\ast}\\
s_{ick}^{\ast}
\end{array}
\right]=
\left[
\begin{array}{c}
\frac{\mu+\gamma}{\mu+\gamma+\omega}\\
0\\
0
\end{array}
\right]
\label{eq:free-disease-equilibrium}
\end{equation}

The endemic equilibrium point $P_{e2}$ is given by Eq. \ref{eq:endemic-equilibrium}

\begin{equation}
P_{e2}=\left[
\begin{array}{c}
s^{\ast}\\
i^{\ast}\\
s_{ick}^{\ast}
\end{array}
\right]=
\left[
\begin{array}{c}
\frac{\beta _{1}+\beta _{2}+\mu }{\alpha(1-\theta)}\\
(\beta_3+\mu)\psi\\
\beta_2\psi
\end{array}
\right],
\label{eq:endemic-equilibrium}
\end{equation}

\noindent where 

\begin{equation}
\psi=\frac{\alpha(1-\theta)(\gamma +\mu) -\left(\beta _{1}+\beta _{2}+\mu \right)\,\left(\gamma +\mu +\omega \right)}
{\alpha(1-\theta )((\beta _{1}+\beta _{2}+\gamma +\mu)(\beta _{3}+\mu) +\beta _{2}\,\gamma)}.
\label{eq:psi}
\end{equation}

The local stability of the equilibrium points are determined with  the Jacobian matrix 
shown in Eq. \ref{eq:JacobianP2}. 

\begin{small}
\begin{equation}
J = \left[
\begin{array}{@{}c@{}c@{}c@{}}
-\alpha(1-\theta)i^{\ast} -(\mu+\gamma+\omega) & -\alpha(1-\theta)s^{\ast} -\gamma & -\gamma \\
\alpha(1-\theta)i^{\ast} & \alpha(1-\theta)s^{\ast} - (\beta_1+\beta_2+\mu) & 0 \\
0 & \beta_2 & -(\beta_3 +\mu) \\
\end{array}
\right]
\end{equation}
\label{eq:JacobianP2}
\end{small}

The disease free equilibrium point $P_{df2}$ is associated with the eigenvalues 
$\lambda_{1}$, $\lambda_{2}$ and $\lambda_{3}$ in Eqs. \ref{eq:eigP2-1}, \ref{eq:eigP2-2} 
and \ref{eq:eigP2-3}, respectively. 

\begin{align}
\lambda_{1} &= -(\gamma +\mu + \omega)
\label{eq:eigP2-1}\\
\lambda_{2} &= -(\beta_1 + \beta_2 + \mu) + \alpha(1-\theta)\dfrac{\mu+\gamma}{\mu+\gamma+\omega}
\label{eq:eigP2-2}\\
\lambda_{3} &= - (\beta_3 + \mu)
\label{eq:eigP2-3}
\end{align}

From Eq. \ref{eq:eigP2-1} and \ref{eq:eigP2-3} it can be seen that 
$\lambda_{1}$ and $\lambda_{3}$ are negative. However, $\lambda_{3}$ 
is negative only if the condition in \ref{eq:eigP2-2-estabilidade} holds true. 

\begin{equation}
\omega>\left(\frac{\alpha(1-\theta)}{\beta_{1}+\beta_{2}+\mu}-1\right)(\mu+\gamma)=
\left(R_{0}-1\right)(\mu+\gamma)
\label{eq:eigP2-2-estabilidade}
\end{equation}

Therefore, considering the foregoing relations, and if condition \ref{eq:eigP2-2-estabilidade} 
holds, the point $P_{e2}$ is locally asymptotically stable. 

Now, regarding the endemic equilibrium point $P_{e2}$, it exists only if the condition expressed 
in \ref{eq:endemic_equil_exist_cond-omega} holds true. 

\begin{equation}
\omega\leq\left(\dfrac{\alpha(1-\theta)}{\beta_1 + \beta_2 + \mu}-1\right)(\mu+\gamma)=
\left(R_{0}-1\right)(\mu+\gamma)
\label{eq:endemic_equil_exist_cond-omega}
\end{equation}

It is interesting to notice that 
the existence condition depends on the vaccination rate $\omega$. Besides,   
for $R_{0}<1$, $P_{e2}$ exists only for negative values of $\omega$, which 
is a case of no physical concern. 

The stability of the endemic equilibrium point  $P_{e2}$ is conditioned to  
parameters combinations, including the social distancing $(\theta)$ and the vaccination 
rate $(\omega)$ parameters. The mathematical reasoning to derive and prove such 
stability conditions is cumbersome, for that reason, numerical analysis is used to infer 
the stability of the point. 

The simulations results presented in Section \ref{sec:SimResults}, 
show that if the disease free equilibrium point $P_{df2}$ is unstable, then the trajectories 
in the vicinity of $P_{e2}$ are stable. These results can be seen from the numerical analysis 
represented by the bifurcation graphs in Section \ref{sec:SimResults}.
 
\section{\SV\ model parametric fitting}
\label{sec:ParamFit}
In this section, the parameters of the  SIRSi model  (Eq. \ref{eq:constant-perc-pop-reduced-model4})  
are numerically determined by parametric fitting using public available data covering a one year 
period of the \cvd\ pandemic, specifically, ranging from March, 2020 to March, 2021. 

The \SV\ model considers the social distancing influence on the \cvd\ infectiousness.  
Recently, a number of works have reported that \cvd\ epidemiological models are strongly sensitive 
to social distancing index \cite{Giordano:2020wq,PREM2020e261,Ferguson2020,BATISTELA2021110388}.  

The social distancing index is time-varying, presenting an oscillatory behavior  
with approximately 7-days period, a feature that impairs the numerical 
parametric adjustment of the \SV\ model. 

In order to avoid numerical problems, the approach is to use a polynomial function 
adjusted to the social distancing data as the social distancing index. The data is 
obtained using privacy-protected mobile phone device location information available 
in SEADE \cite{SEADEISOLAMENTO} for \SP, \Camp\ and \St.

\begin{figure}[!htb]
\centering
\includegraphics[width=\linewidth]{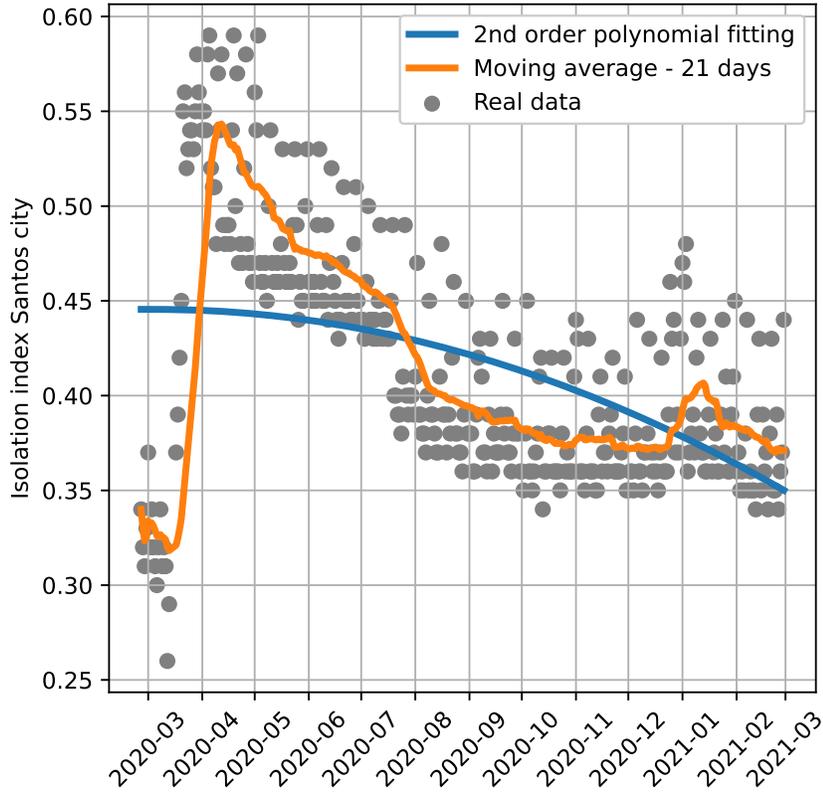}
\caption{Isolation index of \St.}
\label{fig:isol_santos}
\end{figure}

\begin{figure}[!htb]
\centering
\includegraphics[width=\linewidth]{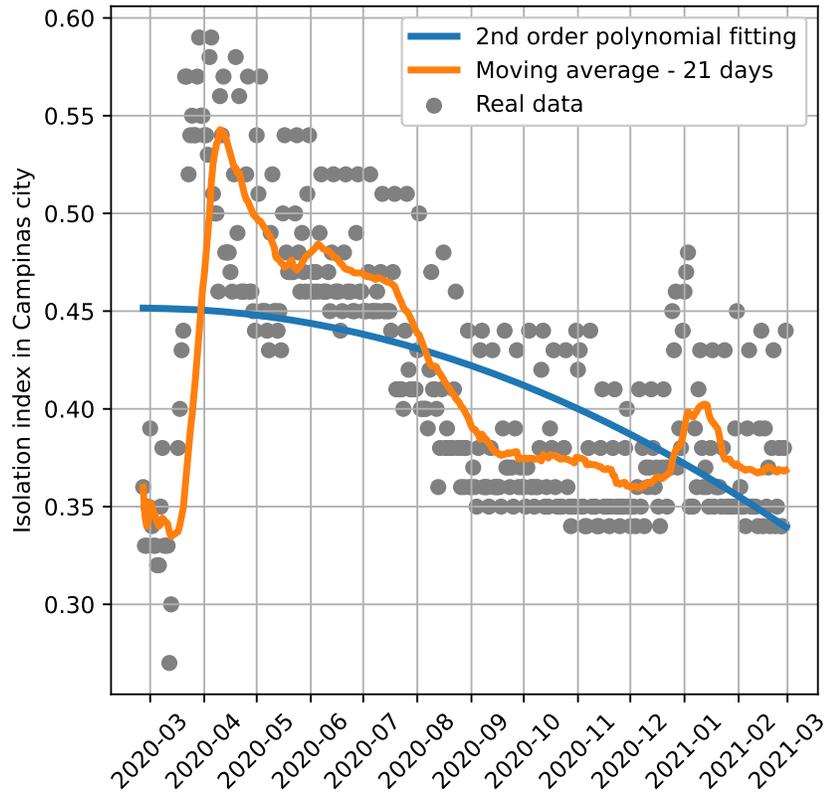}
\caption{Isolation index of \Camp.}
\label{fig:isol_campinas}
\end{figure}

\begin{figure}[!htb]
\centering
\includegraphics[width=\linewidth]{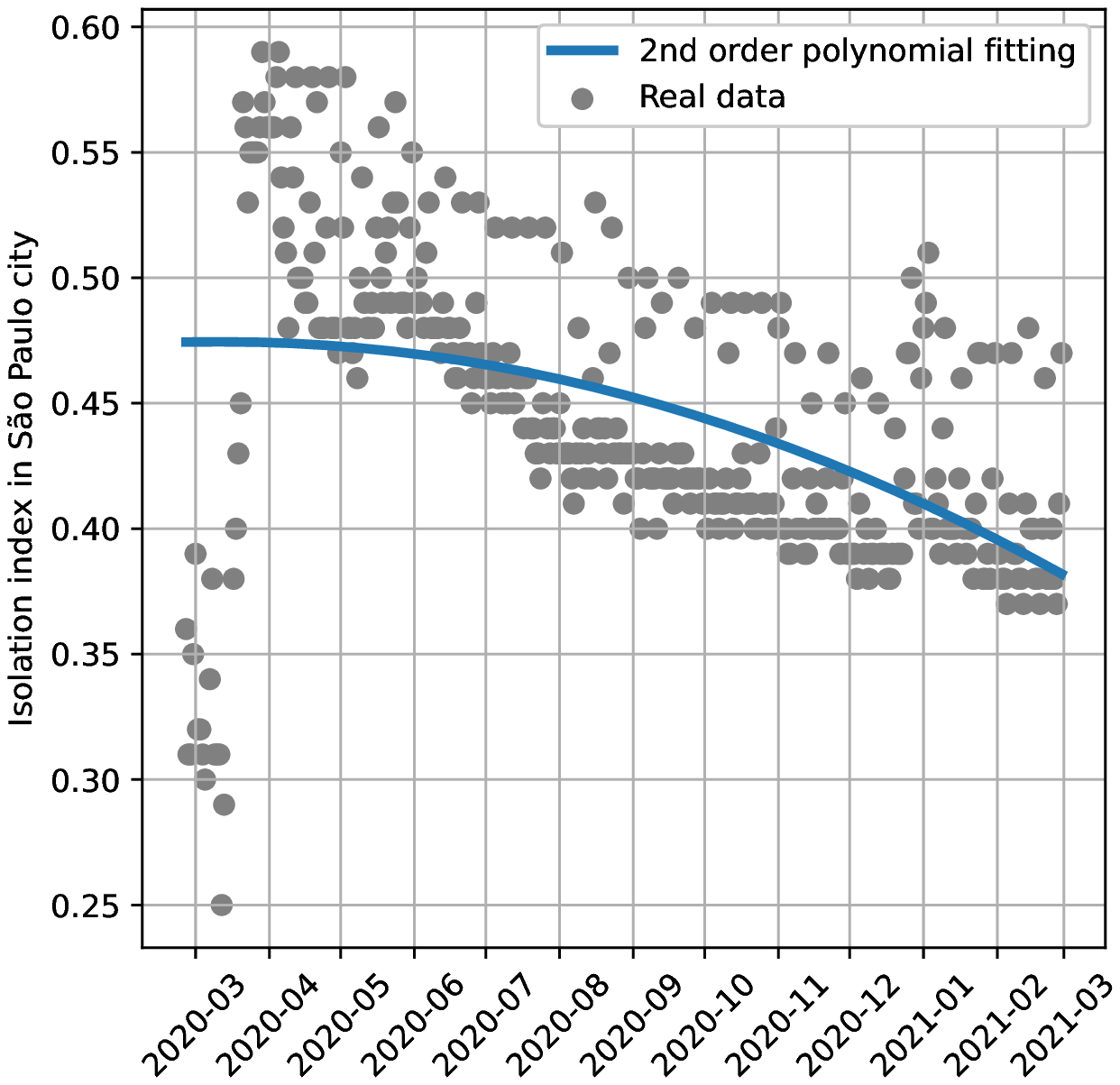}
\caption{Isolation index of \SP.}
\label{fig:isol_saopaulo}
\end{figure}

For the least-squares parameters adjustment raw data was used for \SP. On the other hand, 
for \Camp\ and \St, instead of the raw data, a 21-day moving average filter was used in order 
to smooth the data prior to the least-squares fitting. In spite of that, the filtered data preserves 
both monthly trend and general behavior. The social distancing index polynomial fitting 
results can be seen in Figs. \ref{fig:isol_santos},  \ref{fig:isol_campinas} and  \ref{fig:isol_saopaulo}.

For all the other parameters, the Least-squares Trust Region Reflective (TRR) 
algorithm \cite{doi:10.1137/1.9780898719857,doi:10.1137/0904038},  which is a robust real-time 
optimization method, was used to fit the number of \cvd\ confirmed cases. 

The free parameters in the fitting process are $\left(\gamma, \alpha, \beta_{1}, \beta_{2}, \beta_{3}\right)$ 
(see Fig. \ref{fig:sirsi-model-vacc}), and $s(0)$ which is the normalized amount of susceptible individuals 
prior the pandemic. The other normalized initial conditions are $s_{ick}(0)=0$ and $i(0)=1-s(0)$. 
In addition, for the lower and upper bounds of the parameters, epidemiological and clinical characteristics 
of \cvd\ reported in the literature where considered 
\cite{Souza:2020tk,CAICEDOOCHOA2020316,MUNAYCO2020338,NABI2020110046}

\begin{figure}[!htb]
\centering
\includegraphics[width=\linewidth]{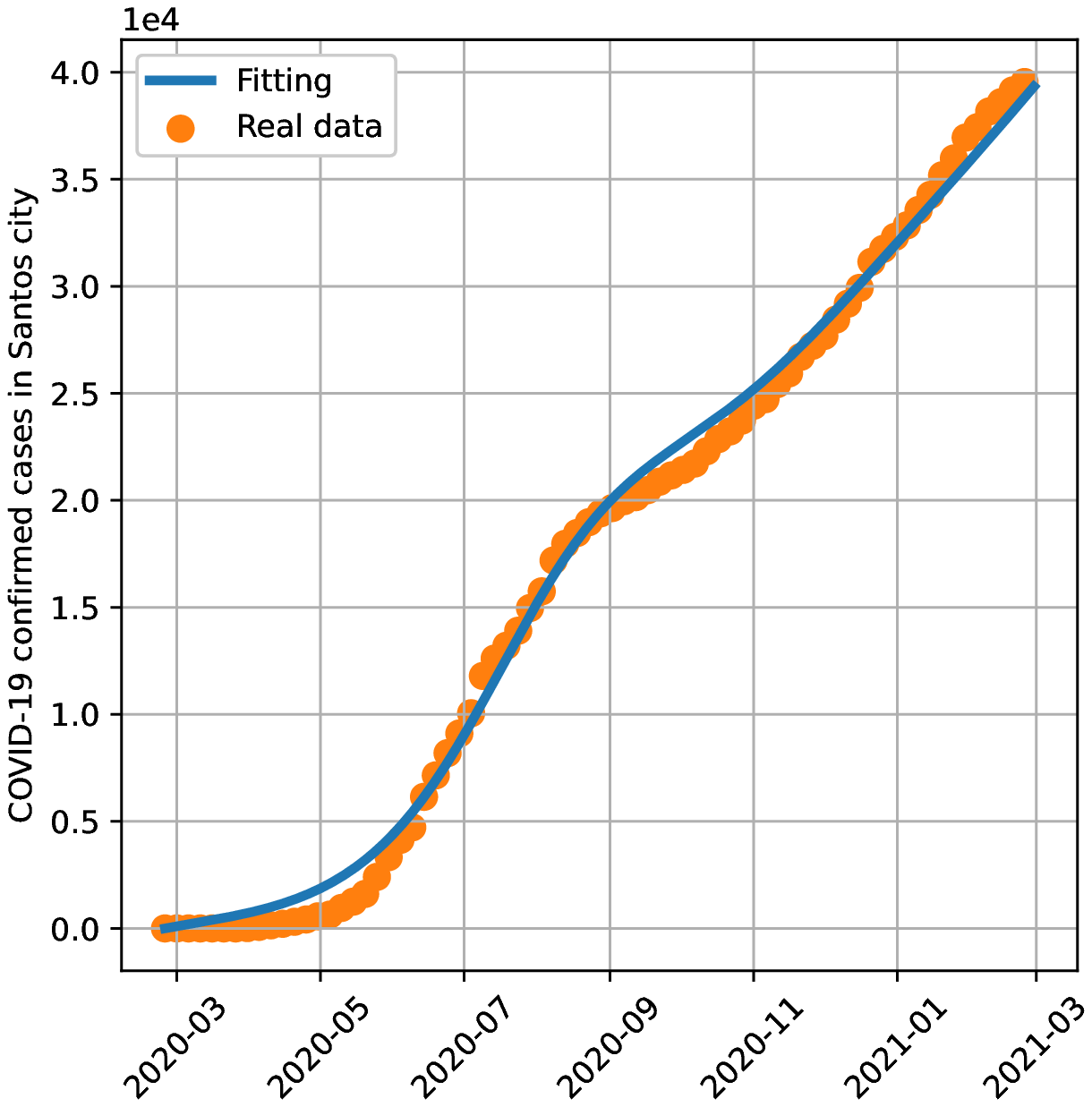}
\caption{\cvd\ infection confirmed cases in \St\ \cite{SEADESocDist} compared to the respective $s_{ick}$ compartment.}
\label{fig:casos_confirmados_santos}
\end{figure}

\begin{figure}[!htb]
\centering
\includegraphics[width=\linewidth]{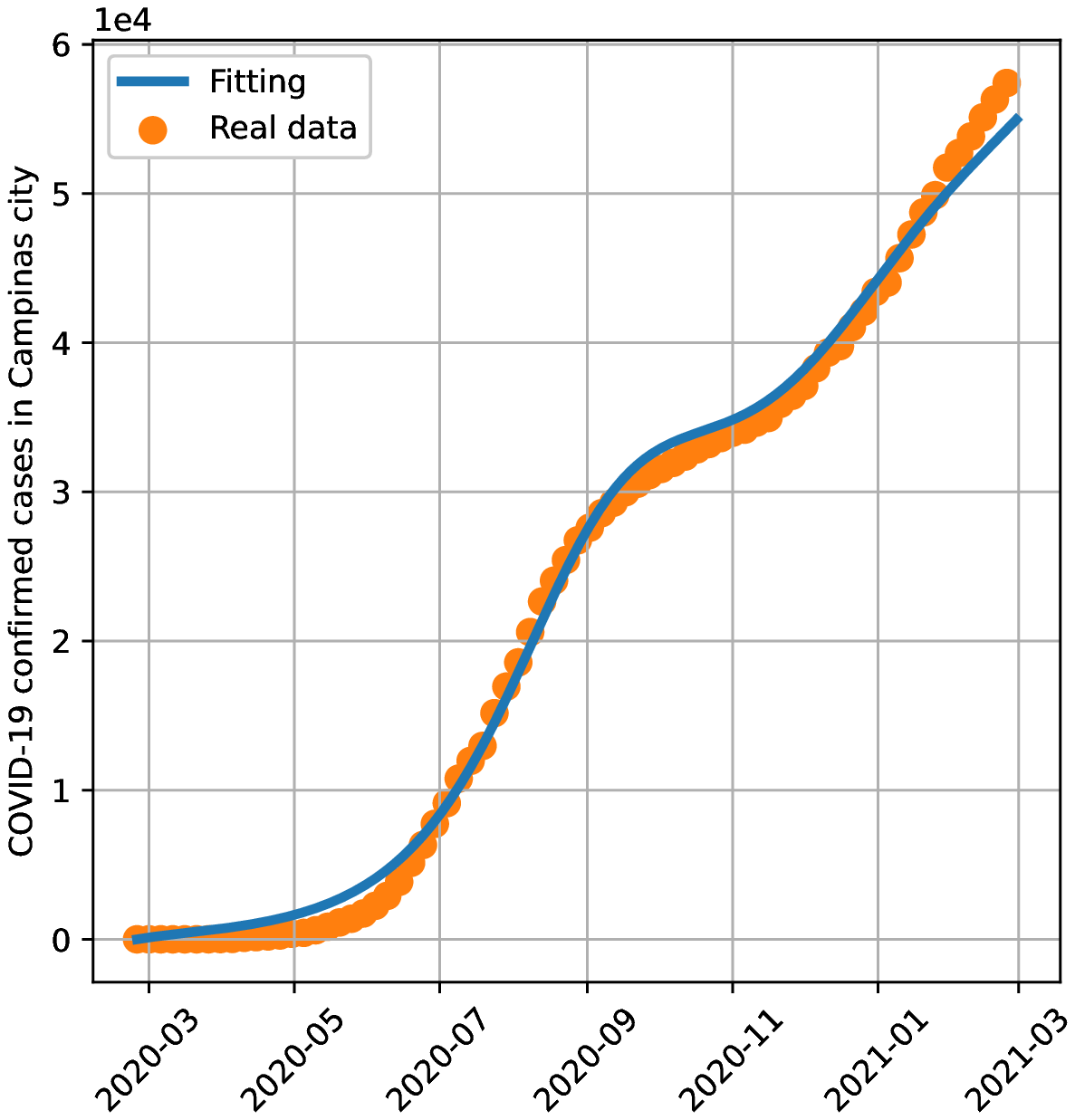}
\caption{\cvd\ infection confirmed cases in \Camp\ \cite{SEADESocDist}  compared to the respective $s_{ick}$ compartment.}
\label{fig:casos_confirmados_campinas}
\end{figure}

\begin{figure}[!htb]
\centering
\includegraphics[width=\linewidth]{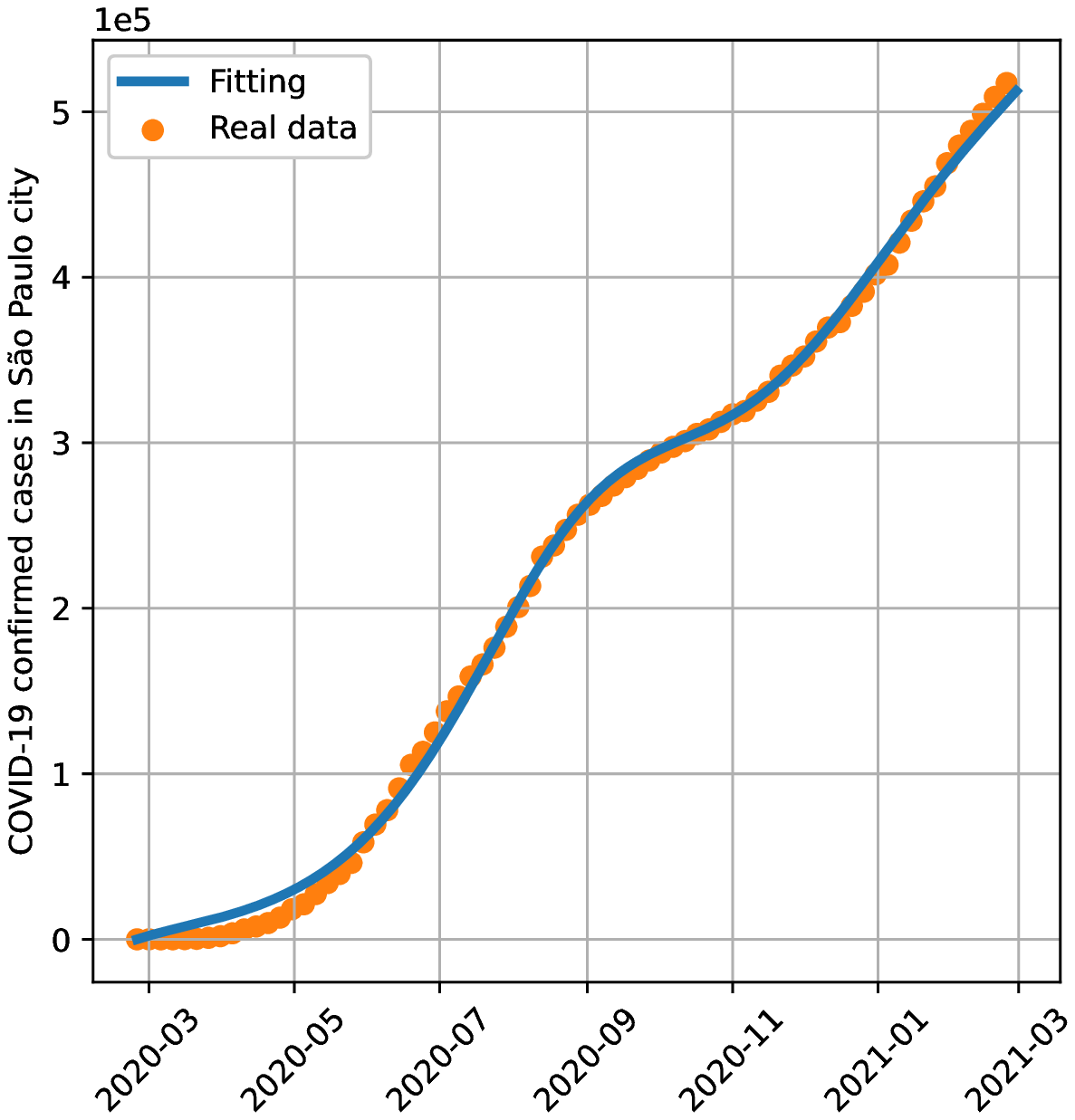}
\caption{\cvd\ infection confirmed cases in \SP\ \cite{SEADESocDist}  compared to the respective $s_{ick}$ compartment.}
\label{fig:casos_confirmados_saopaulo}
\end{figure}

The confirmed cases are considered to be the $s_{ick}(t)$ compartment population in the model of
 Eq. \ref{eq:constant-perc-pop-reduced-model3}. The results of the least-squares TRR fitting 
 comparing the confirmed cases and the respective $s_{ick}$ compartment simulation can 
 be seen in Figs. \ref{fig:casos_confirmados_santos}, \ref{fig:casos_confirmados_campinas} and 
 \ref{fig:casos_confirmados_saopaulo}. 
 
The parameters determined with the least-squares TRR algorithm for \St, \Camp, and \SP, are shown in the 
Tables \ref{tab:fit_santos}, \ref{tab:fit_campinas} and \ref{tab:fit_saopaulo}, respectively. 

\begin{table}[!htb]
\centering
\begin{tabular}{lr}
\hline
Parameter&Value \\
\hline
$\mu$ & 0.000027 \\
$\gamma$ & 0.100000 \\
$\alpha$ & 0.775985 \\
$\theta$ & 0.415355 \\
$\beta_1$ & 0.200000 \\
$\beta_2$ & 0.200000 \\
$\beta_3$ & 0.047847 \\
 $s_0$ & 0.999754 \\
 $i_0$ & 0.000246 \\
\hline
\end{tabular}
\caption{Fitted parameters for \St.}
\label{tab:fit_santos}
\end{table}

\begin{table}[!htb]
\centering
\begin{tabular}{lr}
\hline
Parameter&Value\\
\hline
$\mu$ & 0.000034 \\
$\gamma$ & 0.038255 \\
$\alpha$ & 0.776520 \\
$\theta$ & 0.414454 \\
$\beta_1$ & 0.200000 \\
$\beta_2$ & 0.200000 \\
$\beta_3$ & 0.067000 \\
$s_0$ & 0.999883 \\
$i_0$ & 0.000117 \\
\hline
\end{tabular}
\caption{Fitted parameters for \Camp.}
\label{tab:fit_campinas}
\end{table}

\begin{table}[!htb]
\centering
\begin{tabular}{lr}
\hline
Parameter &    Value \\
\hline
$\mu$ & 0.000036 \\
$\gamma$ & 0.032774 \\
$\alpha$ & 0.811656 \\
$\theta$ & 0.444603 \\
$\beta_1$ & 0.200000 \\
$\beta_2$ & 0.200000 \\
$\beta_3$ & 0.058792 \\
$s_0$ & 0.999800 \\
$i_0$ & 0.000200 \\
\hline
\end{tabular}
\caption{Fitted parameters for \SP.}
\label{tab:fit_saopaulo}
\end{table}

\section{Simulation results}
\label{sec:SimResults}
In this section, the numerical simulations of the \SV\ model (Eq. \ref{eq:constant-perc-pop-reduced-model3}) 
are performed. The parameters used for the \SV\ model simulations are the same 
parameters determined for the SIRSi model in Section \ref{sec:ParamFit}, and shown 
in Tables \ref{tab:fit_santos}, \ref{tab:fit_campinas} and \ref{tab:fit_saopaulo}. 

In the simulations the social distancing index $\theta$ and the vaccination rate $\omega$ 
are kept constant over time. However, they are set up with different values in each 
performed simulation.  

Two different scenarios are considered. Firstly, the vaccination rate is set to $\omega=0$, 
and the numerical simulations are performed for different values of the social distancing index, 
ranging from $\theta=0.0$, i.e., no social distancing intervention, to $\theta=0.7$ which means 
that $70$\% of the population is locked home, keeping social distancing. 

In the next scenario, the social distancing index $\theta$ is set up equal to the average value of the 
measured social distancing data. The values of $\theta$ used for each case can be seen in the 
Tables \ref{tab:fit_santos}, \ref{tab:fit_campinas} and \ref{tab:fit_saopaulo}. In this scenario,  
the simulations are performed considering different values of vaccination rate $\omega$ in order 
to assess the needed vaccination effort to stop, or, at least, to mitigate the spread of the \cvd.  

\subsection{Simulation results for \St}
\label{sec:SimResultsSantos}
\begin{figure}[!htb]
\centering
\includegraphics[width=\linewidth]{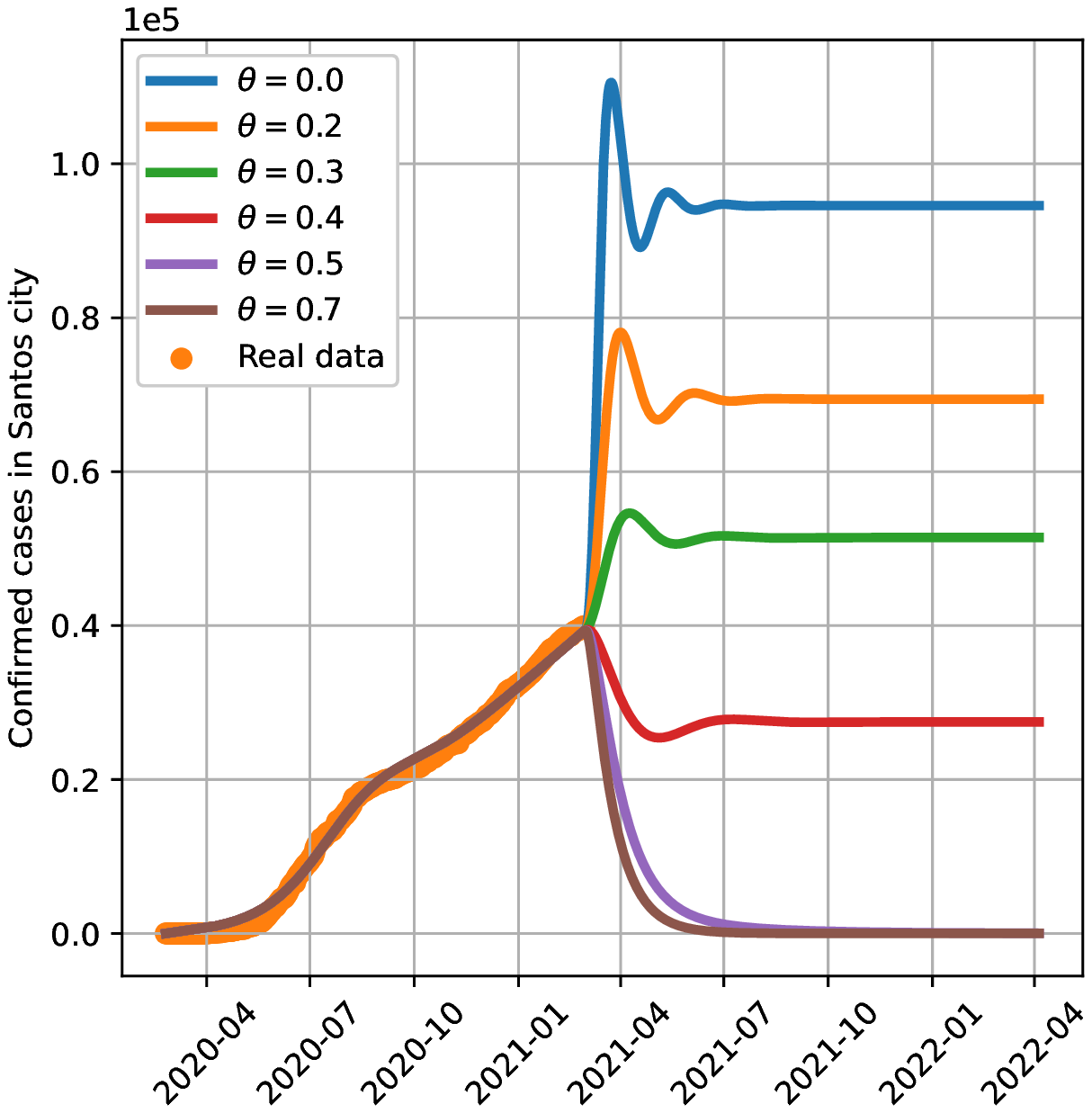}
\caption{Simulation of infection confirmed cases for Santos with $\omega=0$ and different values of social distancing index.}
\label{fig:casos_confirmados_santos_theta_var}
\end{figure}

The simulations in Fig. \ref{fig:casos_confirmados_santos_theta_var} show that without 
vaccination the \cvd\ spread do not decay for a social distancing index less than $\theta=0.5$. 
On the other hand, for $\theta=0$ the number of cases nearly triple, which is a potential 
threath to health care units. With $\theta=0.7$ the infected population rapidly decay until the 
mid July, a behavior that is observed also for $\theta=0.5$, but slower in this last case.
 
It is important to point out that this hypothetical behavior offers an insight to what could 
be the result in such situation, since the social distancing data shows that $\theta$, in fact, 
presents a moving average value behavior combined with a 7-days periodic oscillation.  
  
\begin{figure}[!htb]
\centering
\includegraphics[width=\linewidth]{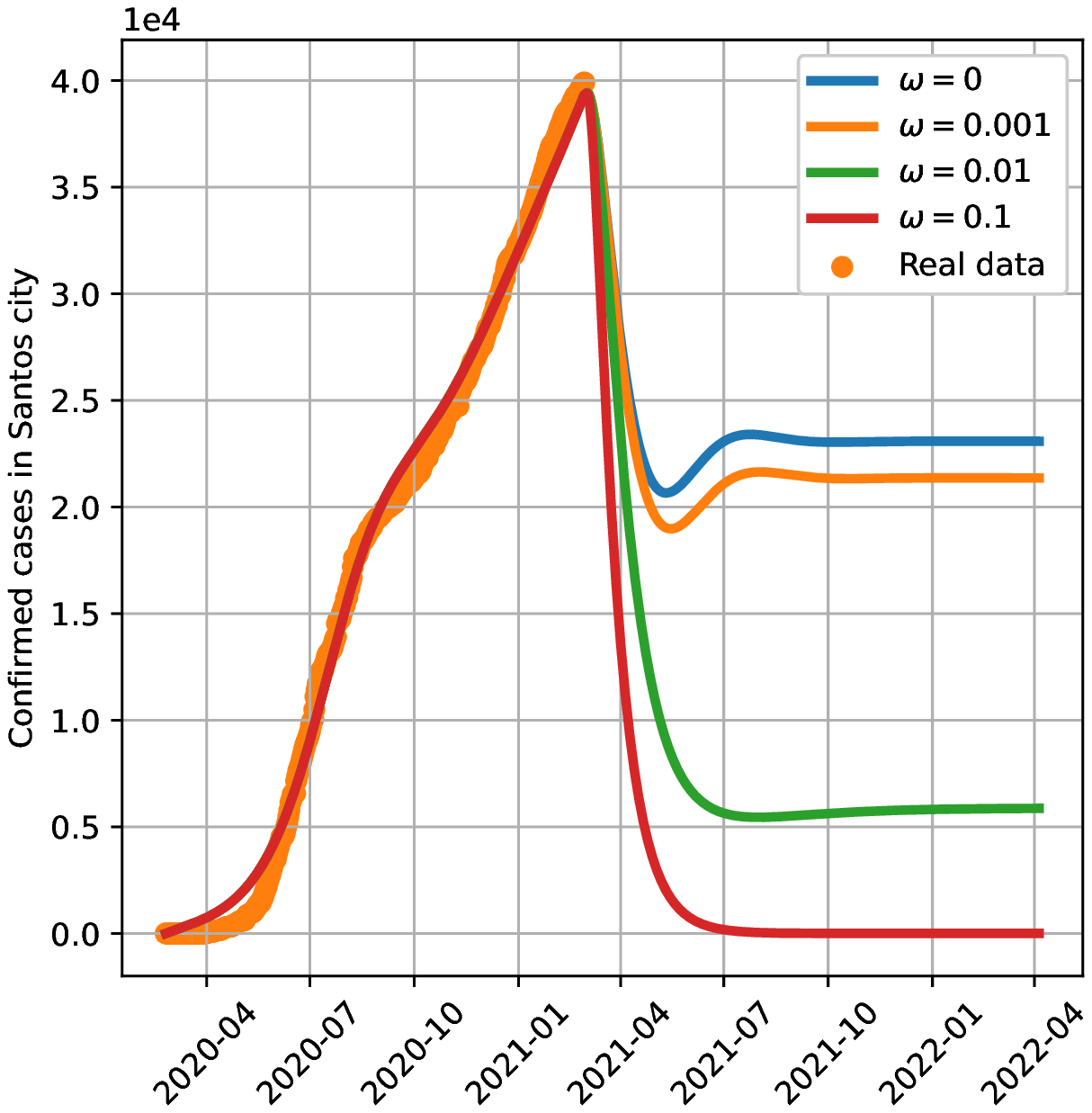}
\caption{Simulation of infection confirmed cases for Santos considering different values of vaccination rate $\omega$. }
\label{fig:casos_confirmados_santos_omega_var}
\end{figure}

The simulations in Fig. \ref{fig:casos_confirmados_santos_omega_var} are performed adopting 
constant social distancing index (see Table \ref{tab:fit_santos}). As it can be seen, the simulations 
results are consistent with the equilibrium points stability and existing conditions in Section \ref{sec:EqPointsSIRSiVaccine}. 

For the vaccination rate $\omega<0.1$ in Fig. \ref{fig:casos_confirmados_santos_omega_var}, the disease 
free equilibrium point $P_{df2}$ is unstable according to condition  \ref{eq:eigP2-2-estabilidade}, 
and the endemic equilibrium point $P_{e2}$ exists according to condition \ref{eq:endemic_equil_exist_cond-omega}. 
The simulation also suggests that $P_{e2}$ is stable due to the steady number of confirmed cases over time. 

For $\omega=0.1$ the situation is the opposite. In this case, the disease free point $P_{df2}$ is stable, 
and the endemic equilibrium point $P_{e2}$ does no exist, according to  condition \ref{eq:endemic_equil_exist_cond-omega}, 
and the number of confirmed cases decay to zero. 

\begin{figure}[!htb]
\centering
\includegraphics[width=\linewidth]{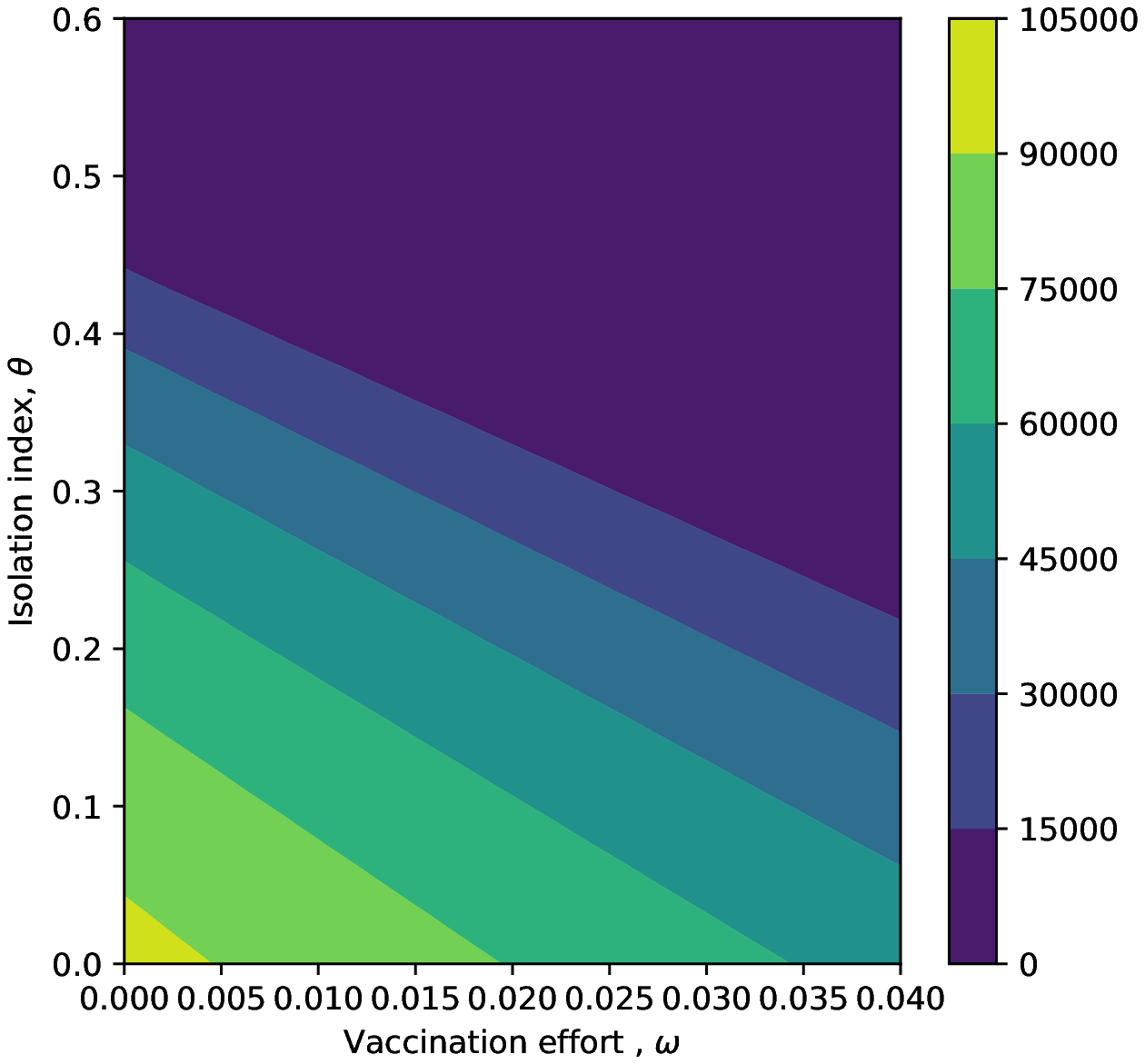}
\caption{Steady state response of the number of confirmed cases for \St.}
\label{fig:omega_x_theta_santos_2}
\end{figure}

In Fig. \ref{fig:omega_x_theta_santos_2} the number of confirmed cases for \St\ as a function of the 
social distancing index $\theta$ and of the vaccination rate $\omega$ is shown. In addition, 
Fig. \ref{fig:omega_x_theta_santos_2} is also a bifurcation diagram, showing the stable 
region (dark-blue) in the parameter space of the disease free equilibrium point. 
The light-blue and green areas indicate the region where the disease free point is unstable. 
The same region indicates the existence of the endemic equilibrium. 
 
In addition, Fig. \ref{fig:omega_x_theta_santos_2} shows the relation between the social distancing 
index $\theta$ and the vaccination rate $\omega$, indicating how to develop a vaccination 
strategy, combined with NPI for the city of \St. 

In this case, with an isolation index $\theta>0.5$ and a small vaccination rate, the number of 
confirmed cases decay to a disease free equilibrium. On the other hand, for $\theta<0.3$, the 
vaccination rate must be much higher, for instance $\omega>0.4$, otherwise the number of 
cases do not decay, remaining constant over time, in a situation compatible with an stable 
endemic equilibrium.

\subsection{Simulation results for \Camp}
\label{sec:SimResultsCampinas}
\begin{figure}[!htb]
\centering
\includegraphics[width=\linewidth]{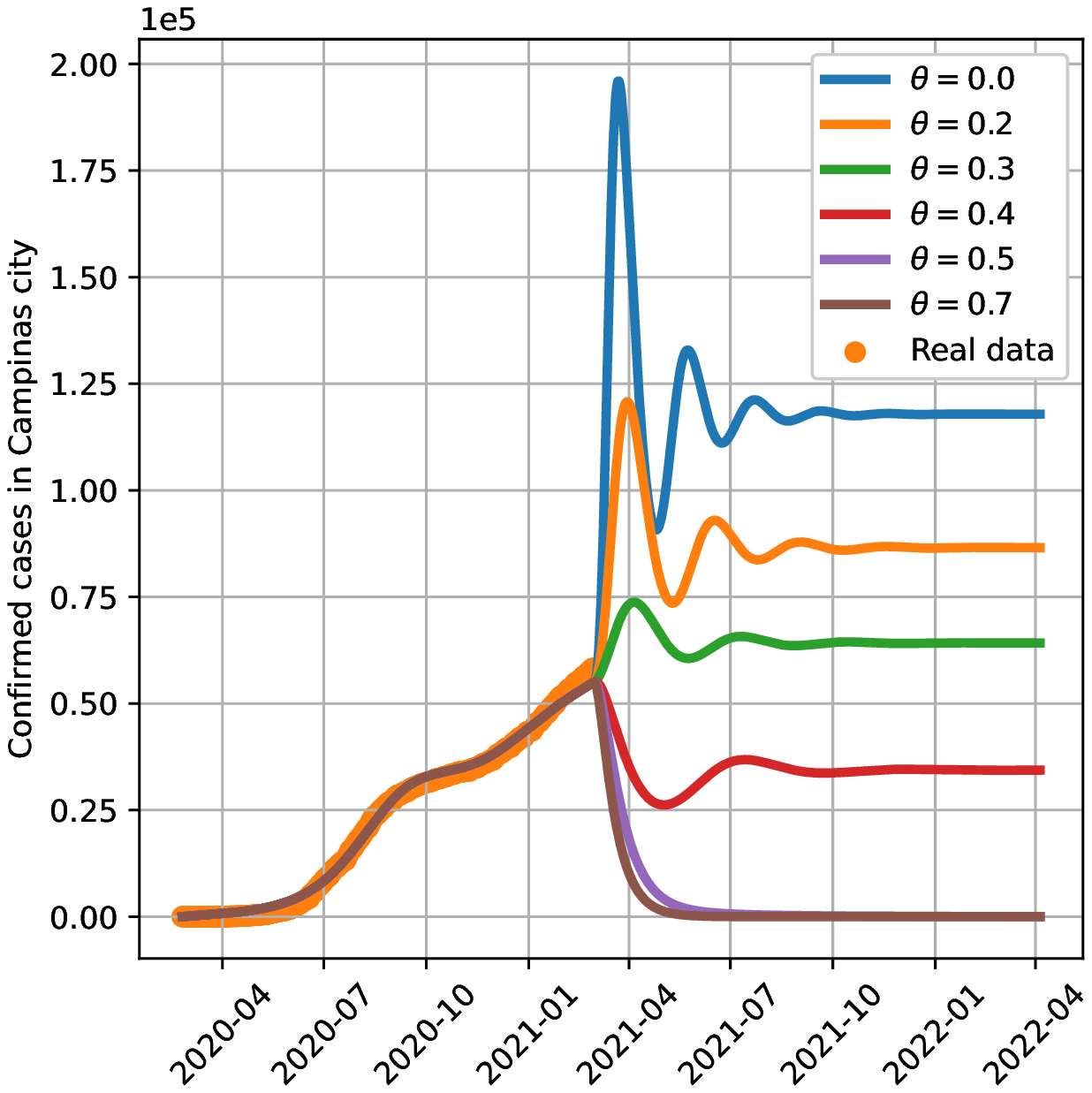}
\caption{Simulation of infection confirmed cases for Campinas with $\omega=0$ and different values of social distancing index.}
\label{fig:casos_confirmados_campinas_theta_var}
\end{figure}

\begin{figure}[!htb]
\centering
\includegraphics[width=\linewidth]{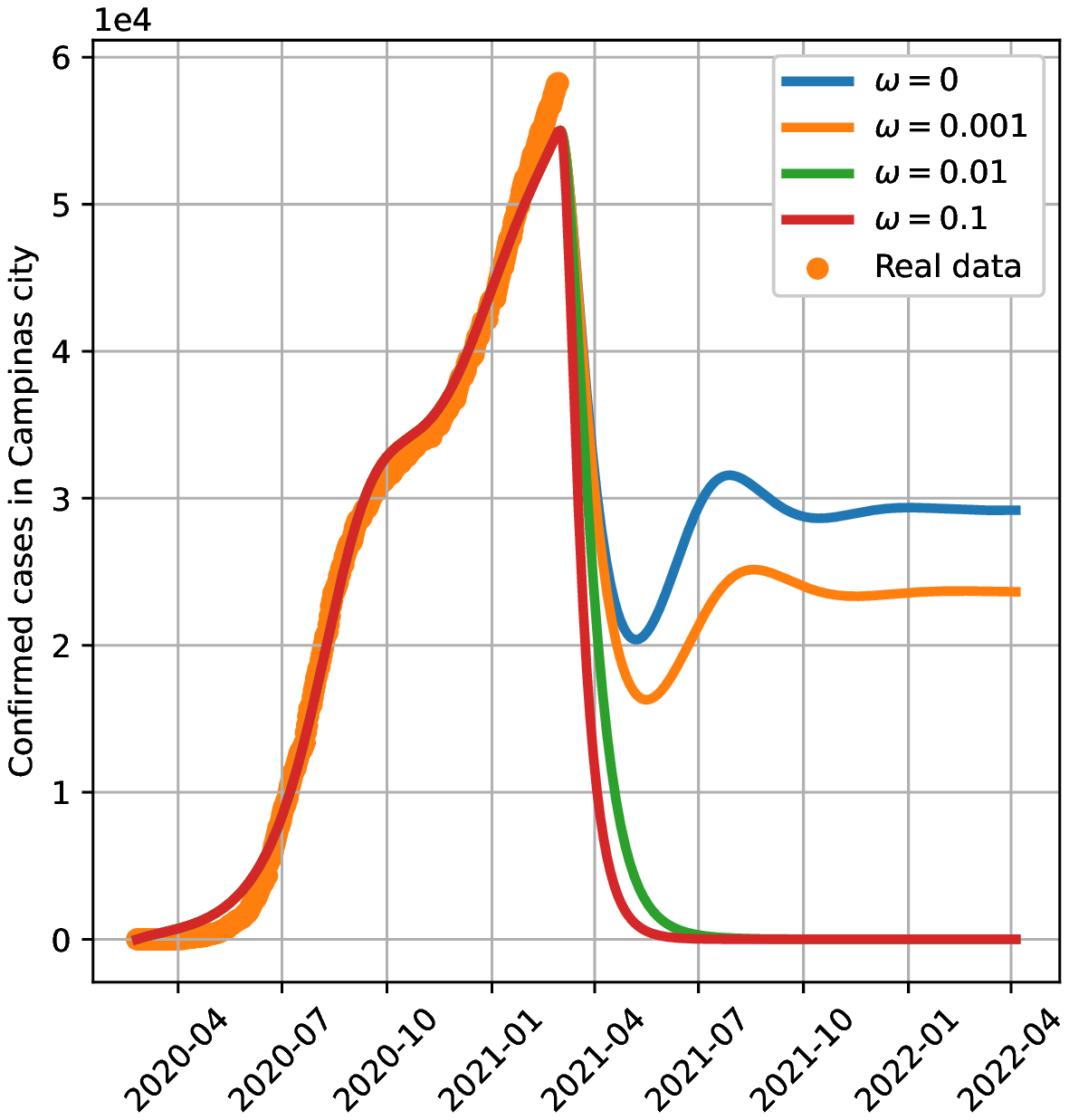}
\caption{Simulation of infection confirmed cases for Campinas considering different values of vaccination rate $\omega$. }
\label{fig:casos_confirmados_campinas_omega_var}
\end{figure}

\begin{figure}[!htb]
\centering
\includegraphics[width=\linewidth]{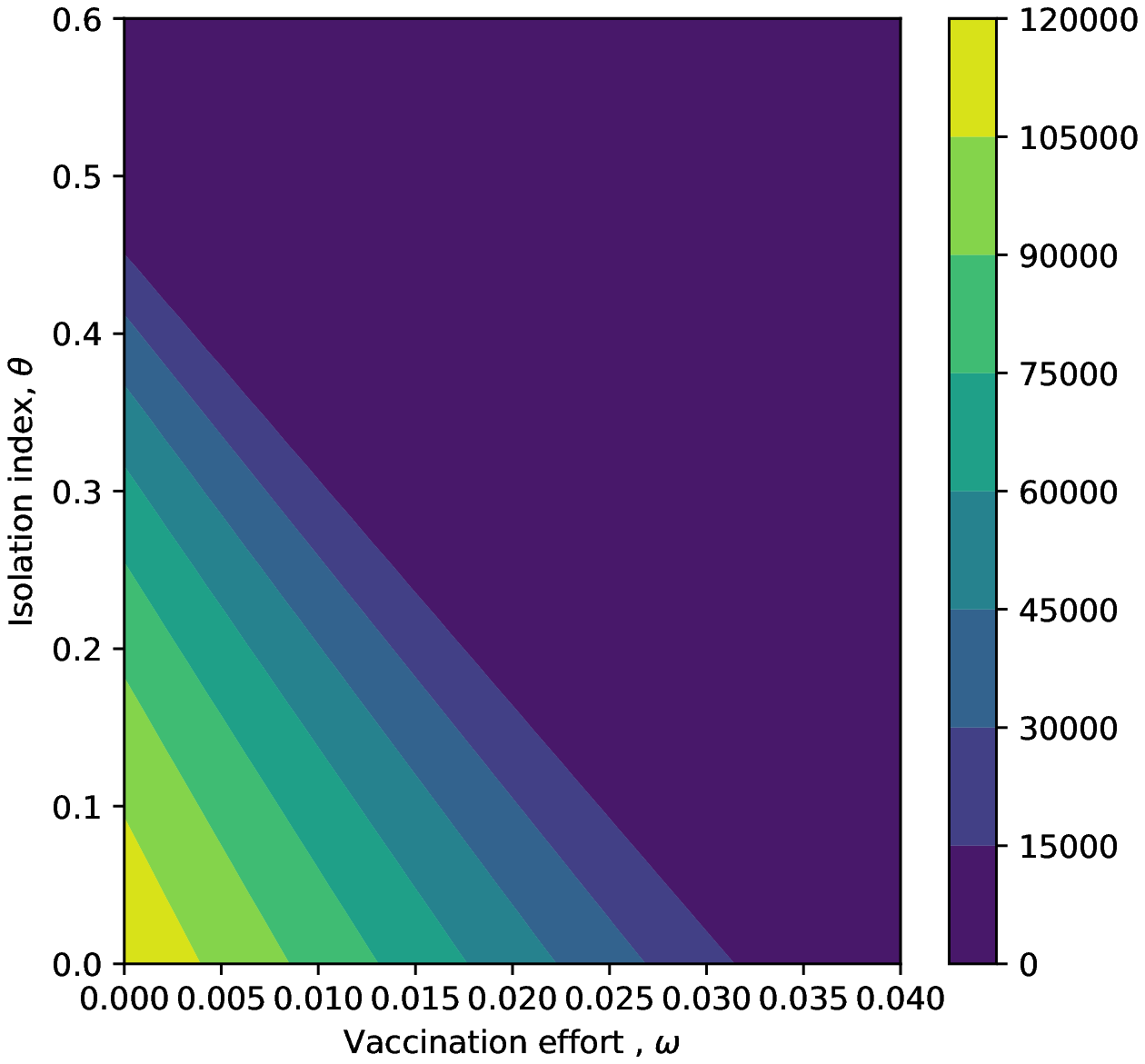}
\caption{Steady state response of the number of confirmed cases for \Camp.}
\label{fig:omega_x_theta_campinas_2}
\end{figure}

The qualitative general behavior of the simulations for \Camp, shown in Figs. \ref{fig:casos_confirmados_campinas_theta_var}, 
\ref{fig:casos_confirmados_campinas_omega_var} and \ref{fig:omega_x_theta_campinas_2},  
is similar to what is presented in the Section \ref{sec:SimResultsSantos}. 

In Fig. \ref{fig:casos_confirmados_campinas_theta_var} the number of confirmed cases is 
shown considering $\omega=0$, and for different values of social distancing index $\theta$. 

As a result, it can be seen that for social distancing index higher than $0.5$, i.e, $\theta>0.5$, 
the disease free equilibrium is stable - condition \ref{eq:eigP2-2-estabilidade} holds true - and 
the number of cases decay to zero over time.  

Differently, for $\theta<0.5$, the disease free equilibrium is unstable. The endemic equilibrium 
exists - condition \ref{eq:endemic_equil_exist_cond-omega} holds true -  and the number of 
cases do not decay over time. It is important to notice that for $\theta<0.3$ the number of confirmed 
cases abruptly increases, nearly quadrupling in a very short period.  This indicates the importance 
of NPI to mitigate the \cvd\ spread and to prevent health care units from collapsing.

The simulations of infected confirmed cases for \Camp, considering a constant social distancing index $\theta$,  
see Table \ref{tab:fit_campinas}, and different values for the vaccination rate $\omega$, are shown in 
Fig. \ref{fig:casos_confirmados_campinas_omega_var}. In this case, a vaccination rate $\omega=0.1$ 
is enough to generate a stable disease free response. For smaller vaccination rates, however, the endemic 
equilibrium exists and the simulations suggests it is stable. 

Comparing Figs. \ref{fig:omega_x_theta_santos_2} and \ref{fig:omega_x_theta_campinas_2}, the vaccination 
rate that is capable of generating stable disease free responses is considerably smaller for \Camp\ than 
for \St. This behavior is probably due to the re-susceptibility feedback gain $\gamma$, since the value found 
for \St\ is nearly $2.6$ times bigger than for \Camp, see Tables \ref{tab:fit_santos} and \ref{tab:fit_campinas}. 

Another relevant information is the slightly bigger social distancing data dispersion for  \St, 
that influences on the infectiousness of the \cvd\ spread, and probably on the 
re-susceptibility feedback gain $\gamma$. 

\subsection{Simulation results for \SP}
\label{sec:SimResultsSaoPaulo}
\begin{figure}[!htb]
\centering
\includegraphics[width=\linewidth]{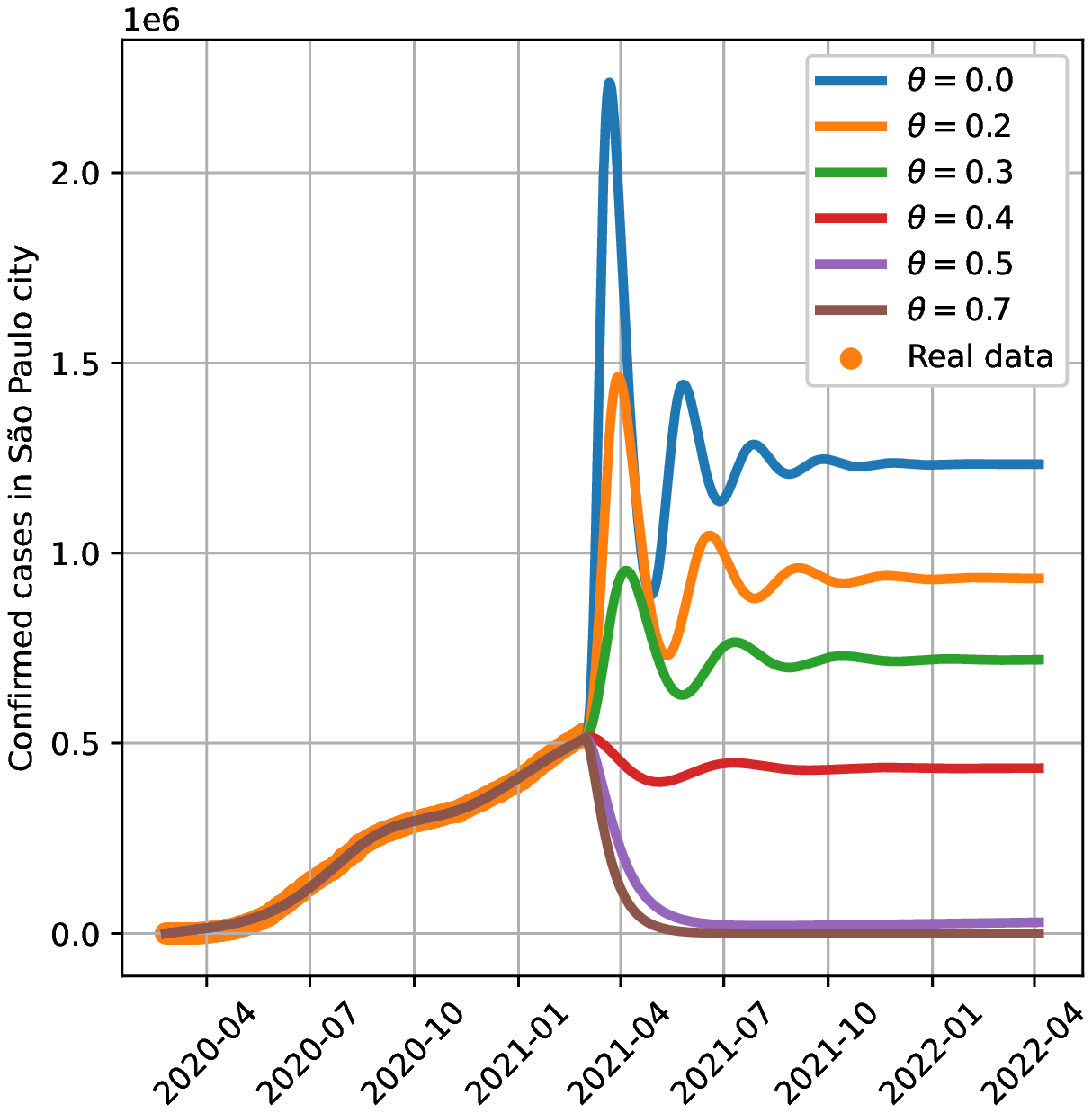}
\caption{Simulation of infection confirmed cases for \SP\ with $\omega=0$ and different values of social distancing index.}
\label{fig:casos_confirmados_saopaulo_theta_var}
\end{figure}

\begin{figure}[!htb]
\centering
\includegraphics[width=\linewidth]{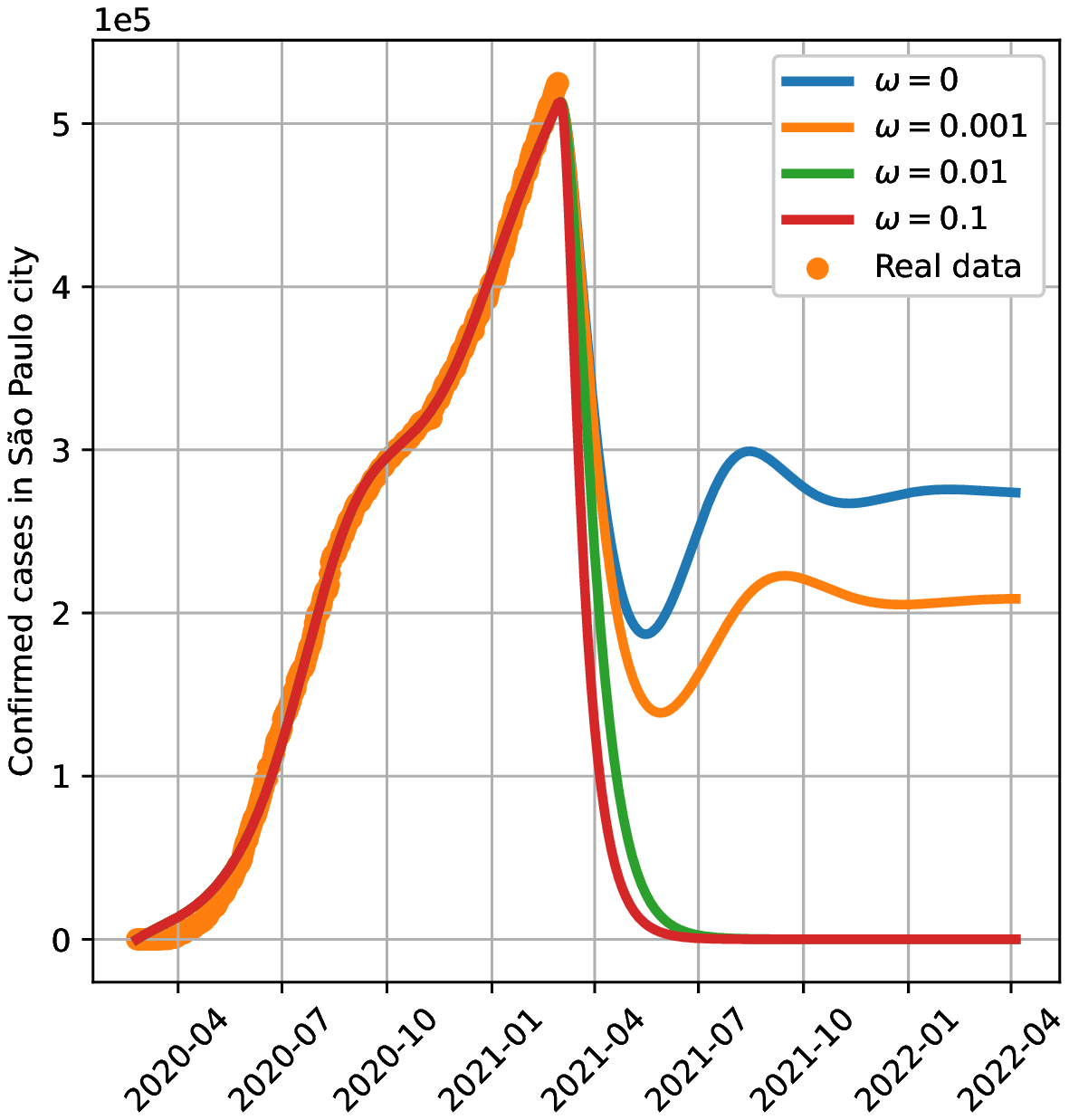}
\caption{Simulation of infection confirmed cases for \SP\ considering different values of vaccination rate $\omega$. }
\label{fig:casos_confirmados_saopaulo_omega_var}
\end{figure}

\begin{figure}[!htb]
\centering
\includegraphics[width=\linewidth]{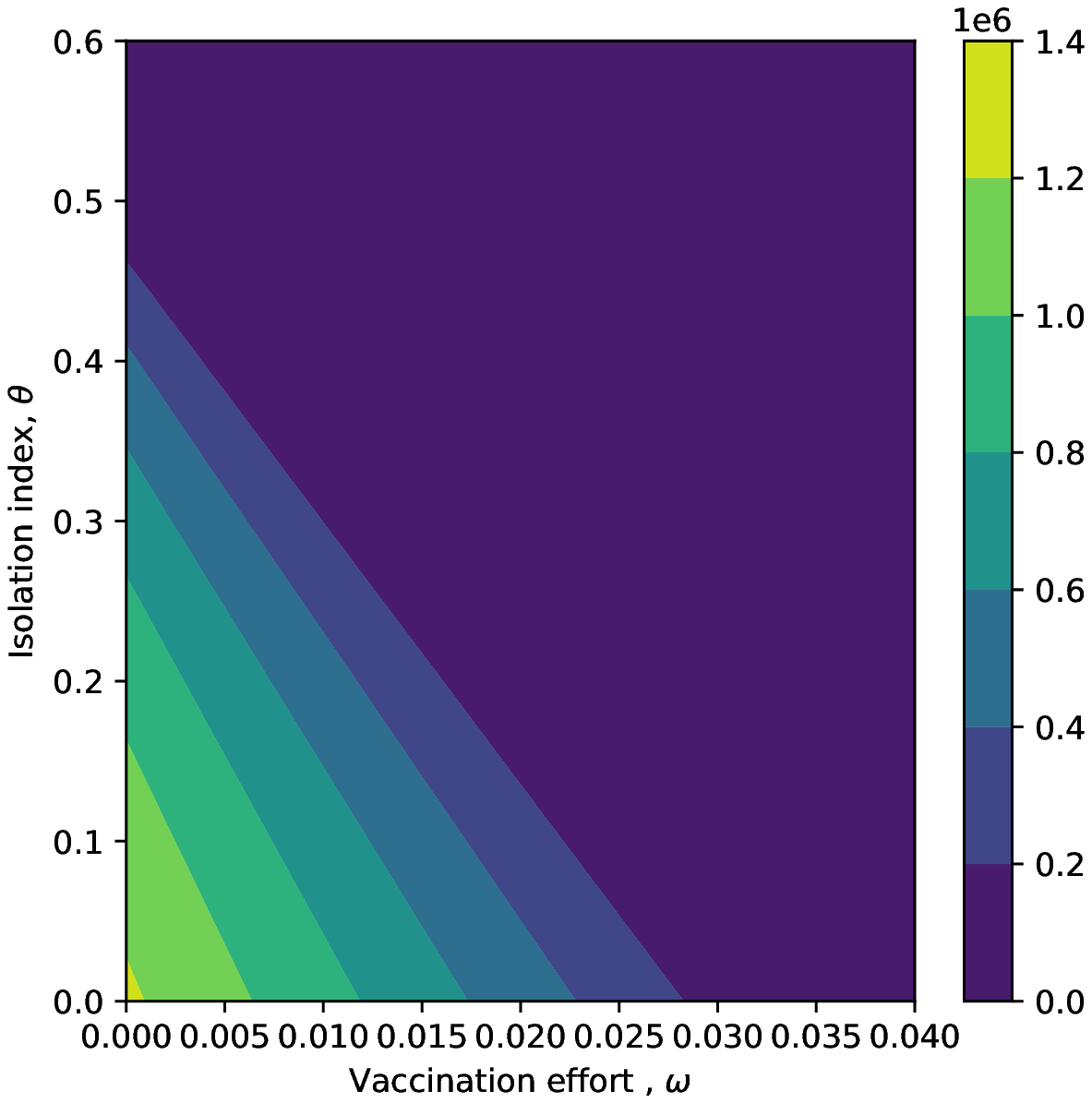}
\caption{Steady state response of the number of confirmed cases for \SP.}
\label{fig:omega_x_theta_saopaulo_2}
\end{figure}

Again, for \SP, the qualitative general behavior of the simulations shown in 
Figs. \ref{fig:casos_confirmados_saopaulo_theta_var}, 
\ref{fig:casos_confirmados_saopaulo_omega_var} and \ref{fig:omega_x_theta_saopaulo_2},  
is similar to the behavior of the simulations presented in the 
Sections \ref{sec:SimResultsSantos} and \ref{sec:SimResultsCampinas}. 

In Fig. \ref{fig:casos_confirmados_saopaulo_theta_var} the simulations of the infection confirmed 
cases for $\omega=0$ and for different values of the social distancing index $\theta$ are shown. 
It can be seen that for $\theta>0.5$ the disease free equilibrium point is stable, i. e., 
condition \ref{eq:eigP2-2-estabilidade} holds true, consequently the number of cases 
decay to zero over time. 

On the other hand, for simulations with $\theta<0.5$ the disease free equilibrium is unstable, and 
since the condition \ref{eq:endemic_equil_exist_cond-omega} holds true,  the endemic equilibrium 
exists. Therefore, the number of cases do not decay over time. 

For $\theta<0.3$ the number of infection confirmed cases increases very fast, threatening the 
health care units of collapsing in a short period of time. 

The simulations of infection confirmed cases shown in Fig. \ref{fig:casos_confirmados_saopaulo_omega_var}, 
indicate that for $\omega\ge0.1$ the disease free equilibrium is stable. On the contrary, for the 
simulated values of $\omega<0.1$, the endemic equilibrium is stable while the disease free equilibrium 
is unstable. As a result, in this case, the number of infection confirmed cases do not decay over time. 

Comparing Figs.   \ref{fig:omega_x_theta_saopaulo_2}, \ref{fig:omega_x_theta_campinas_2} and 
\ref{fig:omega_x_theta_santos_2}, it is clear that the simulation results obtained for \SP\ and 
\Camp\ are very similar concerning the behavior related to the social distancing index $\theta$ 
and the vaccination rate $\omega$. In order to stop the \cvd\ spread in \St, the vaccination 
rate need to be higher. 

The least-square fitting was performed differently for \Camp\ and \St\ due to the dispersion of the 
social distancing index. The influence of this procedure is still not clear, and further investigation 
may be needed in order improve the \SV\ model response.

\section{Conclusions}
\label{sec:Conclusions}
A compartmental model of the \cvd\ pandemic is proposed consisting of four 
compartments,  namely: Susceptible - Infected - Recovered - Sick with vaccination (\SV), considering 
a vaccination strategy and the possibility of recovered patients become susceptible again, that is, 
individuals loose their acquired immunity either after recovery or vaccination. 

The proposed \SV\ model has both disease free and endemic equilibrium points, presenting 
interchangeable stability, that is, the condition that assures stability for the disease free equilibrium 
imply the nonexistence of the endemic equilibrium. On the other hand, the endemic equilibrium 
existence condition imply on the instability of the disease free equilibrium. In addition,  
the existence and stability conditions of expressions \ref{eq:eigP2-2-estabilidade} and 
\ref{eq:endemic_equil_exist_cond-omega} depend on both the reproduction basic number $R_{0}$ 
and on the vaccination rate $\omega$.  

The proposed \SV\ model is least-squares fitted to \cvd\ pandemic publicly available data 
for \St, \Camp\ and \SP. The qualitative general behavior of the simulations is compatible with 
the analytical results. However, comparing the simulations of \St\  with \Camp\ and \SP,  
 a slightly different behavior is found for \St. 
 
The difference is probably due to the re-susceptibility feedback gain $\gamma$, determined by 
the least-squares algorithm, that is $2.6$ times higher for \St, when compared with \Camp\ and 
\SP. This result indicates that for \St,  a successful vaccination strategy imply a 
higher vaccination rate. 

Another result is the influence of the social distancing index $\theta$.  
Small variations on this index can generate abrupt increase in the number of infected confirmed 
\cvd\ cases in a short period of time. This behavior is of special concern due to the possible threat 
to health care units.



\section*{Availability of data and materials}
Data are publicly available  with \cite{SEADEPortalEstatisticas,SEADESocDist}.

\section*{Declaration of competing interest}
There is no conflict of interest between the authors.
\section*{Acknowledge}

JRCP is supported by the Brazilian Research Council (CNPq), grant number: 302883/2018-5.

\end{document}